\documentclass{amsart}

\usepackage{graphicx}

\newtheorem{thm}{Theorem}
\newtheorem{lemma}[thm]{Lemma}
\newtheorem{cor}[thm]{Corollary}

\newcommand{\Nt}{\noindent {\bf Notation: }}
\newcommand{\Pf}{\noindent {\bf Proof: }}
\newcommand{\Rmk}{\noindent {\bf Remark: }}

\begin{document}

\title{Intrinsic Knotting and Linking of Almost Complete Partite
Graphs}        
\author{Thomas W. Mattman, Ryan Ottman, and Matt Rodrigues}        
\address{Department of Mathematics and Statistics,
         California State University, Chico,
         Chico CA 95929-0525, USA}
\email{TMattman@CSUChico.edu}
\address{Department of Mathematics,
	 South Hall Room 6607,
	 University of California,
Santa Barbara, CA 93106}
\email{ryan\_ottman@hotmail.com}
\address{Department of Mathematics, UC Berkeley,
	 970 Evans Hall \#3840,
	 Berkeley, CA 94720-3840}
\email{matthewrodrigues@msn.com}
\subjclass{Primary O5C10, Secondary 57M15}
\keywords{intrinsic linking, embedded graphs, complete partite graphs,
intrinsic knotting}
\thanks{The second and third authors are undergraduate students who 
worked under the supervision of the first author. The second author has
been supported in part by funding from CSU, Chico Research Foundation through the Research
and Creativity Awards program. The first and third author received support from the
MAA's program for Strengthening Underrepresented Minority Mathematics
Achievement (SUMMA) with funding from the NSF and NSA. This research took
place in part during SUMMA REU's at CSU, Chico during the summers of
2003 and 2004.}

\begin{abstract} We classify graphs that are $0$, $1$, or $2$ edges
short of being complete partite graphs with respect to intrinsic
linking and intrinsic knotting. In addition, we classify intrinsic
knotting of graphs on $8$ vertices. For graphs in these families, we
verify a conjecture presented in Adams' {\em The Knot Book}:
If a vertex is removed from an intrinsically knotted graph, one obtains
an intrinsically linked graph.
\end{abstract}

\maketitle

\section{Introduction} 
We say that a graph is intrinsically knotted (respectively, linked)
if every tame embedding of the graph in $\mathbb{R}^3$ contains a non-trivially
knotted cycle (respectively, pair of non-trivially linked cycles).
Robertson, Seymour, and Thomas~\cite{RST}
demonstrated that intrinsic linking is determined by the seven
Petersen graphs. A graph is intrinsically linked if and only if it is or 
contains one of the  seven as a minor. They also showed~\cite{RobSey} that a
similar, finite list of graphs exists for the intrinsic knotting property.
However, determining this list remains difficult. 

It is known~\cite{CG,F1,KS,RMS} that $K_7$ and
$K_{3,3,1,1}$ along with any graph obtained from these two by triangle-Y
exchanges is minor minimal with respect to intrinsic knotting.
Recently Foisy~\cite{F2} has added a new minor minimal graph to the list.
Foisy's example is particularly interesting from our perspective as it
provides a counterexample to the ``unsolved question'' posed in
Adams'~\cite{A} book: {\em Is it true that if $G$ {[}a graph{]} is
intrinsically knotted and any one vertex and the edges coming into it are
removed, the remaining graph is intrinsically linked?} 

While Adams' conjecture is not true in general, it does appear to hold
for a wide array of graphs, particularly graphs that are ``almost" 
complete. We will say a graph is
\textbf{$k$-deficient} if it is a complete or complete partite graph
with $k$ edges removed. In the current
article we verify Adams' conjecture for $0$-, $1$-, and $2$-deficient
graphs. This can be seen as the first few steps in
a project to find a counterexample to Adams' conjecture of minimum
deficiency. However, as Foisy's counterexample is $13$-deficient, there is
a long way to go in this program.

More promising is the search for a counterexample on a minimum number of
vertices. Since $K_7$ is a minor minimal intrinsically knotted graph, it
is the only graph on $7$ or fewer vertices that is intrinsically knotted.
In the current paper, we show that there are twenty intrinsically knotted
graphs on $8$ vertices. These all satisfy Adams' conjecture, so a minimal
counterexample to the conjecture must have between $9$ vertices and the
$13$ of Foisy's graph.

In classifying intrinsic knotting of various families, we have made use of
the known minor minimal examples derived from $K_7$ and $K_{3,3,1,1}$ by
triangle-Y exchanges. In particular, we include here a table of the
$25$ graphs obtained from $K_{3,3,1,1}$ (\cite{KS} includes the table built on
$K_7$). Note that there are no new examples of minor minimal intrinsically knotted 
graphs to be found among $0$-, $1$-, and $2$-deficient graphs and 
graphs on $8$ vertices.

The paper is organized as follows. Following this introduction, Section 2
presents some fundamental lemmas and the table of graphs obtained 
from $K_{3,3,1,1}$ by triangle-Y exchanges. In Sections 3,4, and 5 we
investigate $0$-, $1$-, and $2$-deficient graphs respectively. In each
case we classify the graphs with respect to intrinsic knotting and linking
and demonstrate that they satisfy Adams' conjecture. In Section 6 we
classify intrinsically knotted graphs on $8$ vertices and verify that
they too satisfy Adams' conjecture. We complete the paper by 
formulating the analogue to a question Sachs asked about intrinsic linking.
We show that if $5 \leq n \leq 8$, a graph on $n$ vertices that is not
intrinsically knotted will have at most $5n-15$ edges and ask if 
this is true more generally:

\noindent%
{\bf Question:} Is there a graph on $n$ vertices that is not
intrinsically knotted and has more than $5n-15$ edges?

\section{Lemmas and graphs derived from $K_{3,3,1,1}$}

In this section we present some useful lemmas as
well as a table of graphs derived from $K_{3,3,1,1}$ by triangle-Y 
exchanges. Let us begin with some notation. We will use 
$K_{a_1,a_2,\ldots,a_p}$ to denote a complete partite graph
with $p$ parts containing respectively $a_1, a_2, \ldots, a_p$ vertices.
Permuting the $a_i$ results in the same graph. We will usually 
write the parts in descending order: $a_1 \geq a_2 \geq \ldots \geq a_p$.
The complete graph on $n$ vertices is denoted $K_n$.

\begin{figure}[ht]

\begin{center}
\includegraphics[scale=0.2]{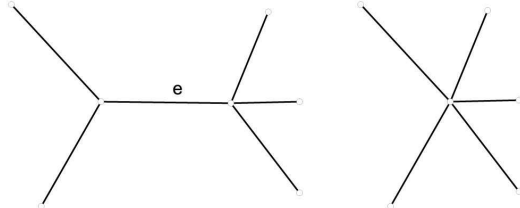}
\caption{An edge contraction of edge e}\label{econt}
\end{center}
\end{figure}

Recall that a  \textbf{minor} of a graph $G$ is the resultant
graph after performing a finite number of vertex or edge deletions and
edge contractions on G. An example of an \textbf{edge contraction} is shown  in
figure \ref{econt} (moving from left to right in the figure). Moving 
from right to left in the figure is known as \textbf{splitting a vertex}.

We say a graph $G$ is \textbf{minor minimal} with respect to a property if 
$G$ exhibits the property but no minor of $G$ does. For example, the 
seven Petersen graphs are the only graphs minor minimal with respect 
to intrinsic linking. An intrinsically linked graph must either be a 
Petersen graph or contain one as a minor.

\begin{figure}[ht]

\begin{center}
\includegraphics[scale=0.74]{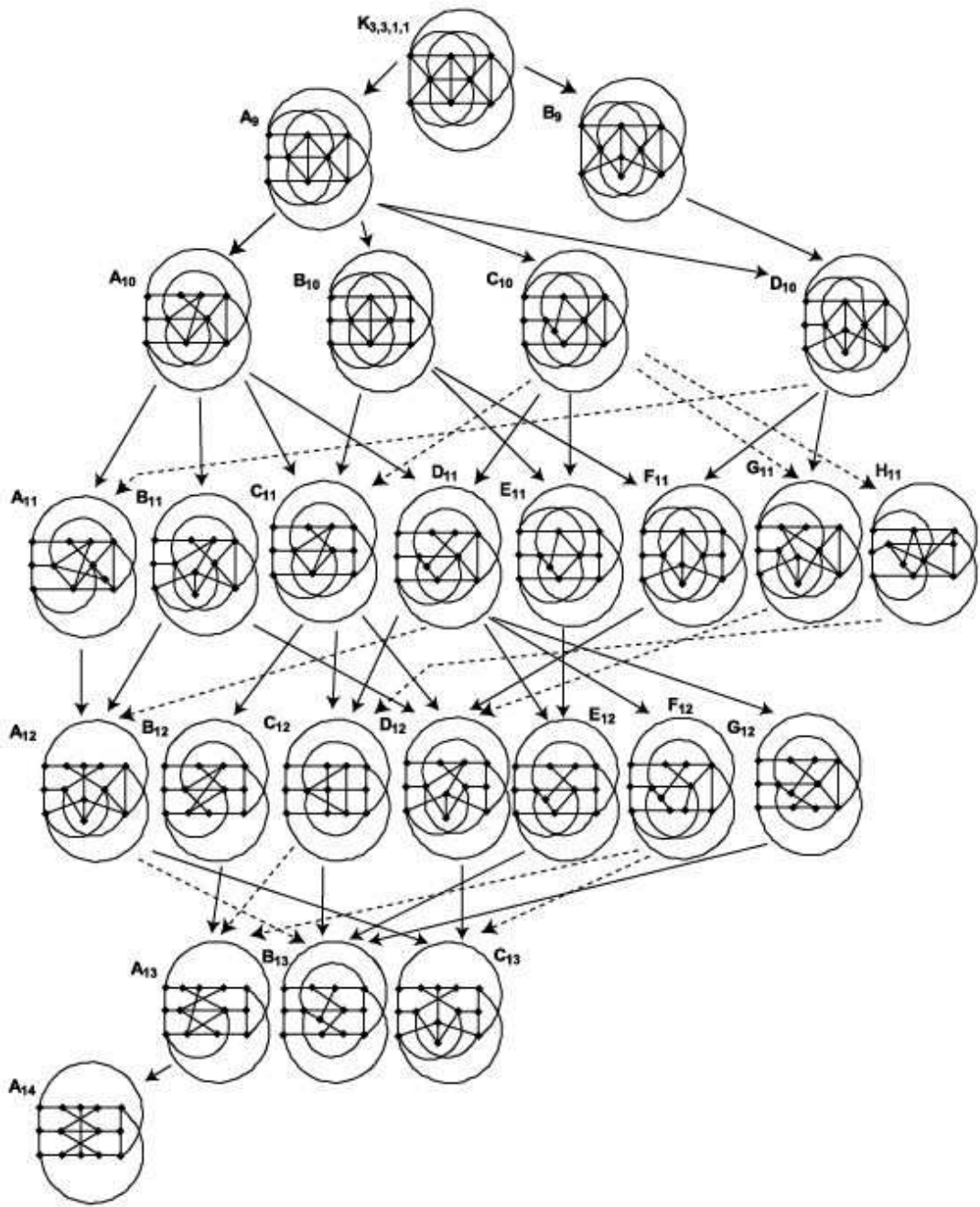}

\caption{Graphs obtained from $K_{3,3,1,1}$ by triangle-Y
exchanges}\label{k3311}
\end{center}
\end{figure}

The analogous list of graphs for
intrinsic knotting is incomplete at this time. 
We do know that the list is finite \cite{RobSey} and 
if $G$ is in the list, so is any graph obtained from
$G$ by triangle-Y exchanges \cite{RMS}. To perform a
\textbf{triangle-Y} exchange on a graph, find 3 vertices that are all 
connected to one another, delete the three edges between them, 
and replace with a
single vertex connected to all 3 vertices of the triangle.
Kohara and Suzuki~\cite{KS} have provided a list of the $13$ graphs
arising from $K_7$ by triangle-Y exchanges. They state that these
graphs are minor-minimal with respect to intrinsic knotting and remark
that there are $25$ graphs that can be constructed from $K_{3,3,1,1}$
by triangle-Y exchanges. We present these graphs in figure~\ref{k3311}. 
Since Foisy~\cite{F1} showed that $K_{3,3,1,1}$
is intrinsically linked it follows \cite{KS} that all $26$ graphs in 
the figure are minor minimal with respect to intrinsic knotting.

The graphs in figure~\ref{k3311} are organized into rows having equal
numbers of vertices. The arrows indicate that a triangle-Y
exchange has occurred. A dashed arrow means that the triangle-Y exchange will
result in a different projection of the graph than that shown in the
figure. Note
that graphs $B_{13}$, $C_{13}$, and $A_{14}$ contain no triangles. 
Most graphs in the figure can be shown to be different from
the others by examining the set of vertex degrees.
The exceptions are the pairs $A_{11}$
and $B_{11}$, $E_{11}$ and $G_{11}$, $A_{12}$ and $B_{12}$, and $A_{13}$
and $B_{13}$. 

We can distinguish the members of these pairs as follows.
$B_{11}$ contains a triangle of degree 5 vertices while
$A_{11}$ does not.
$G_{11}$ contains a degree 3 vertex that is connected only to degree 4
vertices; $E_{11}$ contains no such degree 3 vertex. 
$A_{12}$ contains a 5,5,4 triangle.  The two degree 5 vertices in $B_{12}$
are connected to no common degree 4 vertex so such a triangle is not
possible. 
$A_{13}$ contains a triangle; no triangles exist in $B_{13}$  

These lists of known intrinsically linked or knotted graphs will be 
one of the two main techniques used in our classifications. If we
know that a graph $G$ contains a linked (knotted) minor, then $G$ 
must also be linked (knotted). Conversely, if we can realize $G$
as a minor of an unlinked (unknotted) graph, we deduce that
$G$ must also be unlinked (unknotted). A useful lemma in this 
regard shows how we can combine parts of a $k$-deficient graph.

\begin{lemma}\label{combparts}

$K_{n_1+n_2, n_3, \ldots, n_p}$ is a minor of $K_{n_1,n_2, \ldots, n_p}$.
Similarly, $K_{n_1+n_2, n_3, \ldots, n_p} - k$ edges is a minor of
$K_{n_1,n_2, \ldots, n_p} - k$ edges. 
\end{lemma}

\Pf
Combining the $n_1$  and $n_2$ parts only involves removing edges between
vertices in the $n_1$ part and the $n_2$ part. Recall that with complete
partite graphs, the ordering of the subscripts is not important, so this
lemma implies that we can combine any 2 parts  to get a minor of
the original graph. 

Now, suppose we have a complete partite graph with $k$ edges removed
and we combine two parts. The result would be a complete partite
graph with 1 fewer part and at most $k$ edges removed.  

Furthermore, if there are $m$ edges missing between parts $n_1$ and
$n_2$, we can see that $K_{n_1+n_2, n_3, \ldots, n_k}-(k-m)$
edges is a minor of $K_{n_1,n_2, \ldots, n_k}-k$ edges. 
\qed

The other main technique we will use in our classifications 
is based on a lemma for intrinsic linking due
to Sachs~\cite{S} and an analogous result for intrinsic knotting
proved by Fleming~\cite{Fl}. Let $G + H$ denote the suspension of graphs $G$ and
$H$, i.e., the graph obtained by taking the union of $G$ and $H$ and
adding an edge between each vertex
of $G$ and each vertex of $H$.

\begin{lemma}[\cite{S}]\label{sachs}

The graph $G + K_1$ 
is intrinsically linked if and only if G is non-planar

\end{lemma}

\begin{cor} 
If a vertex is deleted from a graph $H$ and the result is a
planar graph, then $H$ is not intrinsically linked. 
\end{cor}

\Pf 

If the deleted vertex was connected to every other vertex, then by 
lemma \ref{sachs}, $H$ is not intrinsically linked. If the deleted vertex was not
connected to every other vertex, then $H$ is a minor of a graph that is not
intrinsically linked by  lemma \ref{sachs}. \qed

\begin{lemma}[\cite{Fl}]\label{flemming}

The graph $G + K_2$ is intrinsically knotted if and only if
$G$ is non-planar

\end{lemma}

\begin{cor}
If two vertices are deleted from a graph $H$ and the result is a
planar graph, then $H$ is not intrinsically knotted.
\end{cor} 

Note that if we use lemma~\ref{flemming} to argue that a graph 
$H = G+K_2$ is
intrinsically linked, then, $G$ must contain
$K_5$ or $K_{3,3}$ as a minor. It follows that $H$ contains
$K_7$ or $K_{3,3,1,1}$ as a minor. As our only other technique
relies on the family of graphs obtained by triangle-Y exchanges
from $K_7$ and $K_{3,3,1,1}$, our methods cannot possibly add to 
the list of known minor-minimal graphs. In other words,
every intrinsically knotted graph that is $0$-, $1$-, or $2$-deficient
or a graph on $8$ vertices necessarily is or has as a minor one of
the minor minimal graphs obtained from $K_7$ or $K_{3,3,1,1}$ by 
triangle-Y exchanges.

\section{Complete and complete partite graphs}

In this section we classify complete graphs and complete partite graphs 
with respect to intrinsic linking. The classification of these
graphs with respect to intrinsic knotting is due to Fleming~\cite{Fl}.
We use these classifications
to prove that Adams' conjecture holds for this class of graphs.

\begin{thm}

The complete $k$-partite graphs are classified with respect to intrinsic 
linking according to Table \ref{ilcomp}.
\end{thm}

\begin{table}[h]
\begin{center}
\begin{tabular}{|l||c|c|c|c|c|c|} \hline k          & 1 &  2  & 3     &
4       & 5         & $\geq$ 6 \\ \hline  linked     & 6 & 4,4 & 3,3,1 &
2,2,2,2 & 2,2,1,1,1 & All \\ 
           &   &     & 4,2,2 & 3,2,1,1 & 3,1,1,1,1 & \\ \hline not 
linked     & 5 & $n$,3 & 3,2,2 & 2,2,2,1 & 2,1,1,1,1 & None \\
           &   &   & $n$,2,1 & $n$,1,1,1 &           & \\ \hline
\end{tabular}
\end{center}
\caption{Intrinsic linking of complete $k$-partite graphs.}\label{ilcomp}
\end{table}

\Rmk The $n$ in the notation $K_{n,3}$ indicates that the property holds for any
number of vertices in that part, i.e., none of the graphs $K_{1,3}$, $K_{2,3}$,
$K_{3,3}$, $\ldots$ are intrinsically linked.
For each $k$,  the table includes minimal examples of
intrinsically linked complete $k$-partite graphs and maximal graphs which
are not intrinsically linked.
For example, any complete $2$-partite graph which contains
$K_{4,4}$ as a minor is intrinsically linked. On the other hand, any
complete $2$-partite graph which is a minor of a $K_{n,3}$ 
is not linked. Thus,
$K_{l,m}$ is intrinsically linked if and only if $l \geq 4$ and $m \geq 4$.

\Pf

\noindent%
Let us demonstrate that the graphs labeled ``linked" in Table
\ref{ilcomp} are in fact intrinsically linked. 

Conway and Gordon~\cite{CG} proved that
$K_6$ is linked. Any $k$-partite graph with $k \geq 6$ contains  $K_6 =
K_{1,1,1,1,1,1}$ as a minor and is therefore linked. For the remaining
$k$, we appeal to  work of Robertson, Seymour, and Thomas
~\cite{RST} who showed that a graph   is intrinsically linked if and only if
it  is or contains as a minor one of the seven graphs in the Petersen family. In
particular, $K_{4,4}$ with one edge removed and $K_{3,3,1}$ are both in
this family. By combining parts, we see that, for $2 \leq k \leq 5$, one
of these two is a minor of each of the ``linked" graphs in Table 1. 

For each of the ``not linked" examples in Table 1 which involve a part
with a single vertex, the corresponding graph obtained by removing that
vertex is planar. Therefore, by lemma \ref{flemming}, these graphs are
not intrinsically linked. 

Since a cycle requires at least three  vertices, $K_5$ has no disjoint
pair of cycles and is therefore not linked. The remaining ``not linked" 
graphs in Table 1, $K_{n,3}$ and $K_{3,2,2}$, are, respectively, minors
of the unlinked graphs $K_{n,2,1}$ and
$K_{2,2,2,1}$.
\qed

\begin{thm}\label{Aconj} If $G$ is an intrinsically knotted complete
partite graph, and any one vertex and the edges coming into it are
removed, the remaining graph is intrinsically linked.
\end{thm}

\Pf For the reader's convenience, we reproduce Fleming's~\cite{Fl} 
classification of knotted partite graphs as Table \ref{ikcomp} below.

\begin{table}[h]
\begin{center}
\begin{tabular}{|l||c|c|c|c|c|c|c|} \hline k          & 1 &  2  & 3     &
4       & 5         & 6 & $\geq 7$ \\
\hline   knotted    & 7 & 5,5 & 3,3,3 & 3,2,2,2 & 2,2,2,2,1 & 2,2,1,1,1,1
& All \\ 
           &   &     & 4,3,2 & 4,2,2,1 & 3,2,2,1,1 & 3,1,1,1,1,1 &\\
           &   &     & 4,4,1 & 3,3,2,1 & 3,2,1,1,1 &  & \\ 
           &   &     &       & 3,3,1,1 &           &  & \\ \hline

not  knotted    & 6 & 4,4 & 3,3,2 & 2,2,2,2 & 2,2,2,1,1 & 2,1,1,1,1,1 &
None \\
           &   &   & $n$,2,2 & 4,2,1,1 & 2,2,1,1,1 & &\\
           &   &   & $n$,3,1 & 3,2,2,1 & $n$,1,1,1,1 & & \\ 
           &   &   &         & $n$,2,1,1 &         & & \\ 
           &   &   &         & $n$,1,1,1 &         & & \\ \hline
\end{tabular}
\end{center}
\caption{Intrinsic knotting of $k$-partite graphs.}\label{ikcomp}
\end{table}

It suffices to verify the theorem for minimal examples of knotted
$k$-partite graphs, $k = 1,2,3,
\ldots$. 

\begin{description}
\item[$k = 1$] A complete graph $K_n$ is intrinsically knotted iff 
$n \geq 7$. Removing a vertex from $K_7$ produces $K_6$ which is
intrinsically linked.

\item[$k=2$] $K_{5,5}$ is the minor-minimal knotted $2$-partite graph.
(Note that $K_{n,4}$ is not intrinsically knotted for $n\geq 5$. This is
implicit in \cite{Fl}.) Removing a vertex from $K_{5,5}$ yields
$K_{5,4}$ which is intrinsically linked.

\item[$k=3$] There are three minor-minimal knotted $3$-partite graphs:
$K_{3,3,3}$, $K_{4,3,2}$, and $K_{4,4,1}$. Removing a vertex from one of
these results in one of the following linked graphs:
$K_{4,4}$, $K_{3,3,2}$, $K_{4,2,2}$, or $K_{4,3,1}$.

\item[$k=4$] Here the minimal knotted graphs are 
$K_{3,2,2,2}$, $K_{3,3,1,1}$ and $K_{4,2,2,1}$ (The graph
$K_{3,3,2,1}$ listed by Fleming~\cite{Fl} is redundant as it
includes 
$K_{3,3,1,1}$ as a minor.) Removing a vertex from any of these we obtain
one of the linked graphs $K_{3,3,1}$,
$K_{4,2,2}$, $K_{2,2,2,2}$, $K_{3,2,1,1}$, $K_{3,2,2,1}$, or
$K_{4,2,1,1}$.

\item[$k=5$] In this case we must check the knotted graphs
$K_{2,2,2,2,1}$ and $K_{3,2,1,1,1}$. (Fleming's~\cite{Fl}
$K_{3,2,2,1,1}$ is redundant.) Taking a vertex from either of these
results in a linked graph: $K_{2,2,2,2}$, 
$K_{3,2,1,1}$, $K_{2,2,1,1,1}$, $K_{2,2,2,1,1}$,  or $K_{3,1,1,1,1}$.

\item[$k=6$] Here there are two knotted graphs: $K_{2,2,1,1,1,1}$ and
$K_{3,1,1,1,1,1}$. After a vertex is deleted, we're left with one of
these linked graphs: $K_{2,2,1,1,1}$, $K_{3,1,1,1,1}$, or
$K_{2,1,1,1,1,1}$.

\item[$k \geq 7$] All such graphs are intrinsically knotted.  On removing
a vertex, we obtain an $l$-partite graph where $l \geq 6$. All such
graphs are intrinsically linked.
\end{description}
\qed

\section{$1$-deficient graphs}

In this section we classify $1$-deficient graphs with respect to 
intrinsic linking and knotting.  We use the classification to 
prove Adams' conjecture for this family of graphs.

\Nt Often, we will have to talk about a particular vertex or part. We
will refer to parts alphabetically with capital letters and vertices of
those parts with lower case letters. For example, in
$K_{4,3,1}$, we will call the part with 4 vertices part A, the part with
3 vertices part B, and the part with 1 vertex part C. An edge between
parts A and C would be labeled (a,c). 

\subsection{Intrinsic linking}

\begin{thm}

The 1-deficient graphs are classified with respect to intrinsic  linking
according to Table \ref{il1def}.
\end{thm}

\begin{table}[ht]
\begin{center}
\begin{tabular}{|l||c|c|c|c|c|c|c|} \hline k          & 1   &  2      &
3       & 4             & 5               & 6             & $\geq$ 7 \\
\hline  linked     & 7-e & 4,4-e   & 4,3,1-e & 2,2,2,2-e     &
2,2,1,1,1-(b,c) & 2,1,1,1,1,1-e & All \\ 
           &     &         & 3,3,2-e & 3,2,1,1-(b,c) & 3,1,1,1,1-(b,c)
&              & \\  
           &     &         & 4,2,2-e & 4,2,1,1-e     & 4,1,1,1,1-e    
&              & \\  
           &     &         &         & 3,3,1,1-e     & 3,2,1,1,1-e    
&              & \\ 
           &     &         &         & 3,2,2,1-e     & 2,2,2,1,1-e    
&              & \\ \hline

not  linked     & 6-e & $n$,3-e & 3,2,2-e & 2,2,2,1-e     &
2,2,1,1,1-(a,b) & 1,1,1,1,1,1-e & None \\
           &     &       & $n$,2,1-e & $n$,1,1,1-e   & 2,2,1,1,1-(c,d)
&              & \\  
           &     &       & 3,3,1-e   & 3,2,1,1-(a,b) & 3,1,1,1,1-(a,b)
&              & \\  
           &     &       &           & 3,2,1,1-(a,c) & 2,1,1,1,1-e    
&              & \\  
           &     &       &           & 3,2,1,1-(c,d) &                
&              & \\ \hline

\end{tabular}
\end{center}
\caption{Intrinsic Linking of 1 Deficient Graphs.}\label{il1def}
\end{table}

\Rmk

Note that most of these graphs have different types of edges. 
For some graphs, removal of one edge (e.g., (a,b)) will result in a
non-linked graph while the removal of a different edge (e.g., (b,c)) will
result in a linked graph. However, in many cases, removal of any edge
will result in the same categorization. In such cases we simply write
``-e''. For example, no matter which edge we remove from $K_{3,2,2}$, we will
obtain a graph that is not intrinsically linked. 

\Pf

\noindent%
\underline{Linked:}

Let us demonstrate that the graphs labeled ``linked" in Table
\ref{il1def} are in fact intrinsically linked.

$K_7$-e has $K_6$ as a minor. To see this, simply delete a vertex that
the removed edge was attached to. Any 1-deficient graph with 7 or more
parts will contain $K_7$-e (or, equivalently, $K_{1,1,1,1,1,1,1}$-e).

Recall that $K_{4,4}$-e is a Petersen graph, and therefore is
intrinsically linked. Note that by lemma \ref{combparts}, $K_{4,3,1}$-e,
$K_{4,2,2}$-e, $K_{2,2,2,2}$-e, $K_{4,2,1,1}$-e, $K_{3,3,1,1}$-e,
$K_{3,2,2,1}$-e, $K_{4,1,1,1,1}$-e, $K_{3,2,1,1,1}$-e, and
$K_{2,2,2,1,1}$-e all contain $K_{4,4}$-e as a minor, and are therefore
all intrinsically linked.

$K_{3,3,2}$-e has 2 cases: $K_{3,3,2}$-(a,b) and $K_{3,3,2}$-(b,c); both
are intrinsically linked. For $K_{3,3,2}$-(a,b), contract the edge
between vertex b and c to get $K_{3,2,1,1}$. For $K_{3,3,2}$-(b,c),
simply delete vertex c for $K_{3,3,1}$. So in either case, $K_{3,3,2}$-e
is intrinsically linked.

By lemma \ref{combparts}, $K_{3,2,1,1}$-(b,c), $K_{2,2,1,1,1}$-(b,c), and
$K_{3,1,1,1,1}$-(b,c) all contain $K_{3,3,1}$ as a minor and are
therefore intrinsically linked.

There are 2 cases of $K_{2,1,1,1,1,1}$-e: $K_{2,1,1,1,1,1}$-(a,b) and
$K_{2,1,1,1,1,1}$-(b,c). For $K_{2,1,1,1,1,1}$-(a,b), simply delete a to
get $K_6$, and notice that $K_{2,1,1,1,1,1}$-(b,c) is equivalent to
$K_{2,2,1,1,1}$. So both cases are intrinsically linked.

\medskip

\noindent%
\underline{Not Linked:}

$K_6$ (or $K_{1,1,1,1,1,1}$) and $K_{3,3,1}$ are Petersen graphs and
therefore minor minimal with respect to intrinsic linking
\cite{RST}, so $K_6$-e and $K_{3,3,1}$-e are not intrinsically
linked by definition of minor minimal.

$K_{n,3}$-e, $K_{3,2,2}$-e, $K_{n,2,1}$-e, $K_{2,2,2,1}$-e,
$K_{n,1,1,1}$-e, and $K_{2,1,1,1,1}$-e are all not intrinsically linked
before the edge is removed; so, naturally, we can remove an edge and still
have a non-intrinsically linked graph.

$K_{3,2,1,1}$-(c,d) is equivalent to $K_{3,2,2}$ which is not
intrinsically linked.

The remaining graphs all have a vertex that is connected to every other
vertex, and the removal of that vertex results in a planar graph. So, by
lemma \ref{flemming}, none are intrinsically linked.
\qed

\subsection{Intrinsic knotting}

\begin{table}[h]
\begin{center}
\begin{tabular}{|l||c|c|c|c|c|c|c|c|} \hline k          & 1   &  2    &  
3     &    4    &     5     &      6      &  7                &  $\geq 8$
\\
\hline   knotted    & 8-e & 5,5-e & 3,3,3-e & 3,2,2,2-e & 2,2,2,2,1-e &
2,2,1,1,1,1-(b,c) & 2,1,1,1,1,1,1-e & All \\ 
           &     &       & 4,3,2-e & 4,2,2,1-e & 3,2,1,1,1-(b,c) &
3,1,1,1,1,1-(b,c) &        & \\ 
           &     &       & 4,4,1-e & 3,3,2,1-e & 4,2,1,1,1-e &
3,2,1,1,1,1-e &                   & \\ 
           &     &     &           & 4,3,1,1-e & 3,3,1,1,1-e &
2,2,2,1,1,1-e &                   & \\ 
           &     &     &           &           & 3,2,2,1,1-e &
4,1,1,1,1,1-e &                   & \\ \hline


not  knotted    & 7-e & n,4-e & 3,3,2-e &  3,3,1,1-e & 3,2,1,1,1-(c,d) &
2,2,1,1,1,1-(a,b) & 1,1,1,1,1,1,1-e & None \\ 
           &     &     & $n$,2,2-e &  2,2,2,2-e & 3,2,1,1,1-(a,b) &
2,2,1,1,1,1-(c,d) &     & \\ 
           &     &     & $n$,3,1-e &  3,2,2,1-e & 3,2,1,1,1-(a,c) &
3,1,1,1,1,1-(a,b) &                   & \\ 
           &     &     &         &  $n$,2,1,1-e &  2,2,2,1,1-e  &
2,1,1,1,1,1-e &                   & \\ 
           &     &     &       &            &  $n$,1,1,1,1-e    & 
&                   & \\ \hline

\end{tabular}
\end{center}
\caption{Intrinsic knotting of 1 deficient graphs.}\label{ik1def}
\end{table}

\begin{thm}

The 1-deficient graphs are classified with respect to intrinsic  knotting
according to Table \ref{ik1def}.
\end{thm}

\Pf

\noindent%
\underline{Knotted:}

Let us demonstrate that the graphs labeled ``knotted" in Table
\ref{ik1def} are in fact intrinsically knotted.

$K_8$-e is has $K_7$ as a minor; to see this, simply delete a vertex that
the removed edge was connected to. All 1-deficient graphs with 8 or more
parts will contain $K_{1,1,1,1,1,1,1,1}$-e = $K_8$-e.

$K_{5,5}-e$ is intrinsically knotted as proved in \cite{Sh}. It
is also helpful to see that it is an expansion of $H_9$ from
\cite{KS}, a fact pointed out in \cite{F1}.

\begin{figure}[ht]
\begin{center}
\includegraphics[scale=0.35]{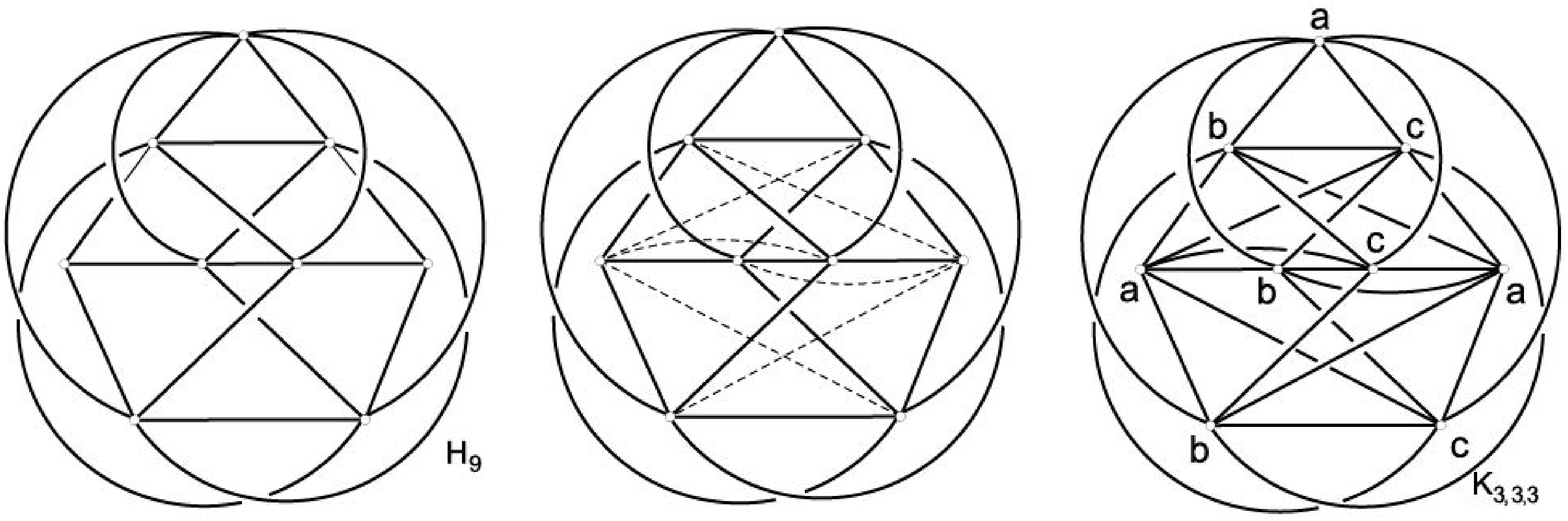}

\caption{$H_9$ is a Minor of $K_{3,3,3}$}\label{h9tok333}
\end{center}
\end{figure}

$K_{3,3,3}-e$ is intrinsically knotted. As shown in figure \ref{h9tok333},
several edges are added to $H_9$ to get $K_{3,3,3}$. $H_9$ was shown to
be intrinsically knotted in \cite{KS}. So, if we simply add 1
fewer of those edges, we arrive at $K_{3,3,3}-e$ and our graph is still
intrinsically knotted. Notice that any two edges in $K_{3,3,3}$ are
equivalent to one another.

By lemma \ref{combparts}, $K_{3,3,2,1}$-e, $K_{3,3,1,1,1}$-e,
$K_{3,2,2,1,1}$-e, $K_{3,2,1,1,1,1}$-e, and $K_{2,2,2,1,1,1}$-e all have
$K_{3,3,3}$-e as a minor, and are therefore all intrinsically knotted.

$K_{4,3,2}-e$ is intrinsically knotted, there are 3 cases.

\begin{figure}[ht]
\begin{center}
\includegraphics[scale=0.35]{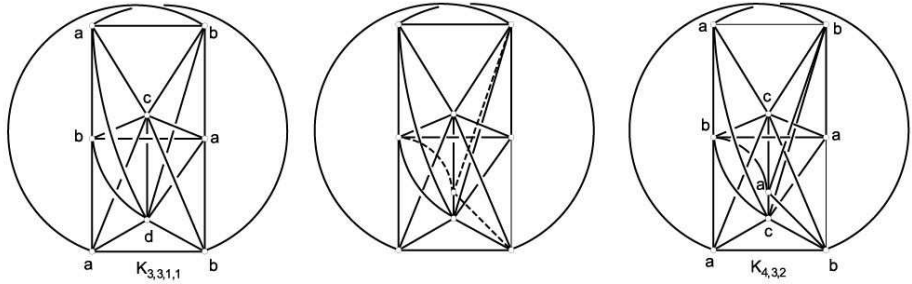}
\caption{$K_{3,3,1,1}$ is a Minor of $K_{4,3,2}$}\label{k3311tok432}
\end{center}

\end{figure}

Case 1: $K_{4,3,2}-(a,b)$. $K_{4,3,2}$ is intrinsically
knotted as it  contains 
$K_{3,3,1,1}$ as a minor (see figure~\ref{k3311tok432}). Notice that if we add
one fewer edge, we get $K_{4,3,2}-(a,b)$.

\begin{figure}[ht]
\begin{center}
\includegraphics[scale=0.35]{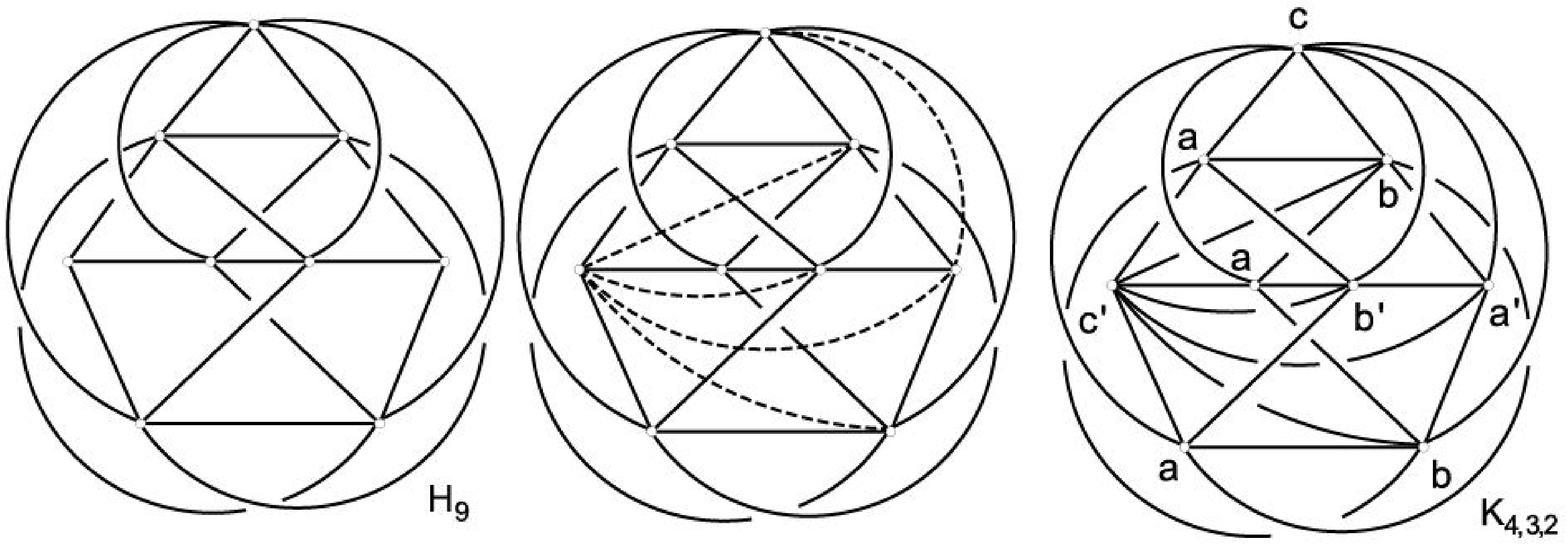}

\caption{$H_9$ is a Minor of $K_{4,3,2}$}\label{h9tok432}
\end{center}
\end{figure}

Case 2: $K_{4,3,2}-(b,c)$. In figure \ref{h9tok432}, $K_{4,3,2}$ is shown to
contain $H_9$ from
\cite{KS}. If we leave out edge
(b$^\prime$,c$^\prime$), we see that $K_{4,3,2}-(b,c)$ contains $H_9$ as
a minor; therefore it is intrinsically knotted.

Case 3: $K_{4,3,2}-(a,c)$. The same as case 2, except we leave out edge
(a$^\prime$,c$^\prime$).

By lemma \ref{combparts}, $K_{3,2,2,2}$-e, $K_{4,2,2,1}$-e,
$K_{4,3,1,1}$-e, $K_{2,2,2,2,1}$-e, $K_{4,2,1,1,1}$-e, and
$K_{4,1,1,1,1,1}$-e all contain $K_{4,3,2}$-e as a minor, and are
therefore intrinsically knotted.

$K_{4,4,1}-e$ is intrinsically knotted; there are 2 cases:
$K_{4,4,1}-(a,b)$ and $K_{4,4,1}-(a,c)$.

\begin{figure}[ht]
\begin{center}
\includegraphics[scale=0.35]{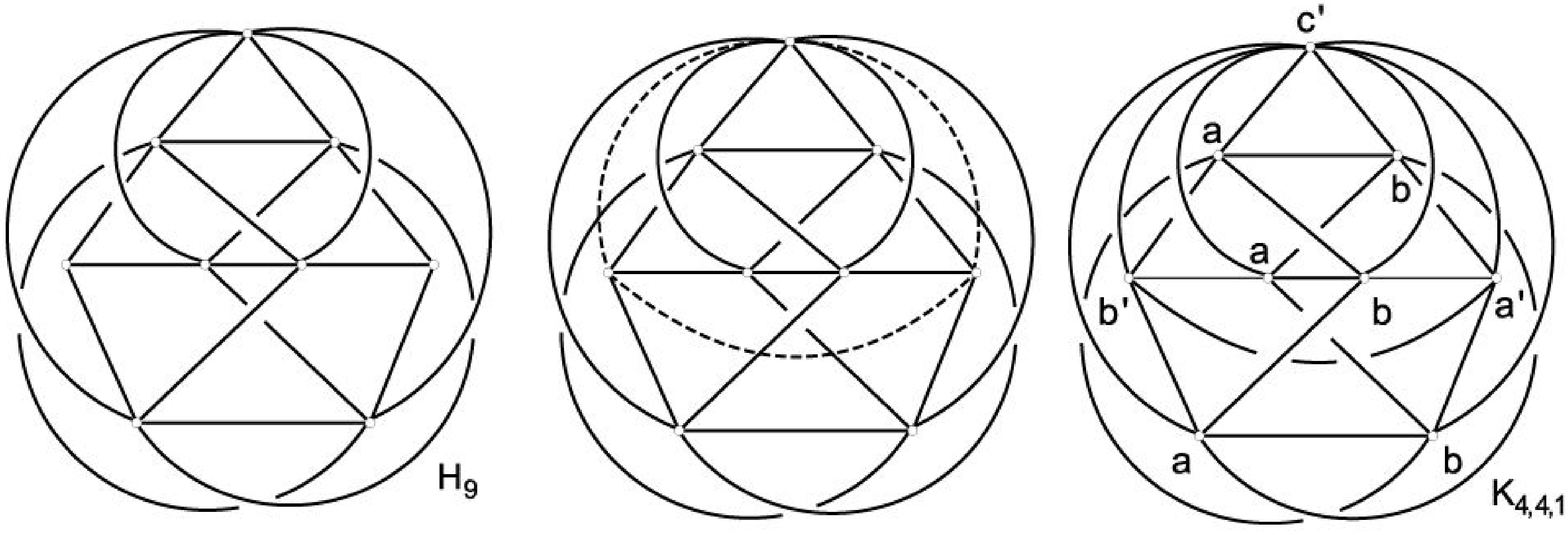}

\caption{$H_9$ is a Minor of $K_{4,4,1}$}\label{h9tok441}
\end{center}
\end{figure}

As shown in figure \ref{h9tok441}, $K_{4,4,1}$ contains $H_9$ as a minor.
If we refrain from adding either (a$^\prime$,b$^\prime$) or
(a$^\prime$,c$^\prime$), we can see that both cases of $K_{4,4,1}$-e have
$H_9$ as a minor, and are therefore intrinsically knotted.

By (the proof of) lemma \ref{combparts}, if we combine the parts of a missing edge
of a 1-deficient graph, we will get a complete graph. If we combine parts B
and C of $K_{3,2,1,1,1}$-(b,c), $K_{2,2,1,1,1,1}$-(b,c), or
$K_{3,1,1,1,1,1}$-(b,c), we will get $K_{3,3,1,1}$ or $K_{3,2,1,1,1}$. So
these three graphs are intrinsically knotted.

$K_{2,1,1,1,1,1,1}$-e has 2 cases: $K_{2,1,1,1,1,1,1}$-(a,b) and
$K_{2,1,1,1,1,1,1}$-(b,c). For $K_{2,1,1,1,1,1,1}$-(a,b), simply delete
vertex a to get $K_7$. $K_{2,1,1,1,1,1,1}$-(b,c) is equivalent to
$K_{2,2,1,1,1,1}$. So both cases are intrinsically knotted.

\medskip

\noindent%
\underline{Not Knotted:}

$K_7$ was shown to be minor minimally intrinsically knotted in
\cite{CG} and $K_{3,3,1,1}$ was shown to be minor minimally
intrinsically knotted in \cite{F1}, so $K_7$-e (or
$K_{1,1,1,1,1,1,1}$-e) and $K_{3,3,1,1}$-e are not intrinsically knotted.

In Fleming's paper \cite{Fl}, he showed that $K_{4,4}$ is not
intrinsically knotted.  In fact, his proof really showed that $K_{n,4}$ is not
intrinsically knotted, so, naturally, $K_{n,4}$-e is not intrinsically
knotted either. Similarly, $K_{3,3,2}$-e, $K_{n,2,2}$-e, $K_{n,3,1}$-e,
$K_{2,2,2,2}$-e, $K_{3,2,2,1}$-e, $K_{n,2,1,1}$-e, $K_{2,2,2,1,1}$-e,
$K_{n,1,1,1,1}$-e, and $K_{2,1,1,1,1,1}$-e are not intrinsically knotted.

$K_{3,2,1,1,1}$-(c,d) is equivalent to $K_{3,2,2,1}$.

The remaining 5 graphs each have 2 vertices connected to every other
vertex. If those 2
vertices are deleted, the resulting graph is planar, so by lemma \ref{flemming},
they are not intrinsically knotted.
\qed

\subsection{Proof of Adams' conjecture for $1$-deficient graphs}

\begin{thm}\label{Aconj1def} If $G$ is a 1-deficient graph, and any one
vertex and the edges coming into it are removed, the remaining graph is
intrinsically linked.
\end{thm}

\Pf

It suffices to verify the theorem for minimal examples of knotted
1-deficient graphs. 

\begin{description}
\item[$k = 1$] $K_8$-e is intrinsically knotted. If we remove a vertex,
we get either $K_7$ or $K_7$-e, both of which are intrinsically linked.

\item[$k=2$] $K_{5,5}$-e is intrinsically knotted. If we remove a vertex,
we get either $K_{5,4}$-e or $K_{5,4}$, both of which are intrinsically
linked.

\item[$k=3$] $K_{3,3,3}$-e, $K_{4,3,2}$-e, and $K_{4,4,1}$-e are
intrinsically knotted. If we remove a vertex from any of these graphs we
will get $K_{3,3,2}$-e, $K_{4,2,2}$-e, $K_{4,3,1}$-e, $K_{4,4}$-e, or one
of these without the removed edge. In any case, the result is
intrinsically linked.

\item[$k=4$] $K_{3,2,2,2}$-e, $K_{4,2,2,1}$-e, $K_{3,3,2,1}$-e, and
$K_{4,3,1,1}$-e are intrinsically knotted. If we remove a vertex from any of these
graphs we will get $K_{4,3,1}$-e, $K_{3,3,2}$-e, $K_{4,2,2}$-e,
$K_{2,2,2,2}$-e, $K_{4,2,1,1}$-e, $K_{3,3,1,1}$-e, $K_{3,2,2,1}$-e, or
one of these without the edge removed. In all cases, the result is
intrinsically linked.

\item[$k=5$] $K_{2,2,2,2,1}$-e, $K_{3,2,1,1,1}$-(b,c), $K_{4,2,1,1,1}$-e,
$K_{3,3,1,1,1}$-e, $K_{3,2,2,1,1}$-e are intrinsically knotted. If we
remove a vertex, we will get $K_{2,2,2,2}$-e, $K_{3,2,1,1}$-(b,c),
$K_{4,2,1,1}$-e, $K_{3,3,1,1}$-e, $K_{3,2,2,1}$-e, $K_{2,2,1,1,1}$-(b,c),
$K_{3,1,1,1,1}$-(b,c), $K_{4,1,1,1,1}$-e, $K_{3,2,1,1,1}$-e,
$K_{2,2,2,1,1}$-e, or one of these without the removed edge. In all
cases, the result is intrinsically linked.

\item[$k=6$] $K_{2,2,1,1,1,1}$-(b,c), $K_{3,1,1,1,1,1}$-(b,c),
$K_{3,2,1,1,1,1}$-e, $K_{2,2,2,1,1,1}$-e, and $K_{4,1,1,1,1,1}$-e are
intrinsically knotted. If we remove a vertex, we will get
$K_{2,2,1,1,1}$-(b,c), $K_{3,1,1,1,1}$-(b,c), $K_{4,1,1,1,1}$-e,
$K_{3,2,1,1,1}$-e, $K_{2,2,2,1,1}$-e, $K_{2,1,1,1,1,1}$-(b,c),
$K_{3,1,1,1,1,1}$-e, $K_{2,2,1,1,1,1}$-e, or one of these without the
removed edge. In all cases, the result is intrinsically linked.

\item[$k=7$] $K_{2,1,1,1,1,1,1}$-e is intrinsically knotted. If we remove
a vertex, we will get $K_{2,1,1,1,1,1}$-e, $K_{1,1,1,1,1,1,1}$-e, or one
of these without the removed edge. All will be intrinsically linked.

\item[$k\geq 8$] All 1-deficient graphs with 8 or more parts are
intrinsically knotted. If we delete a vertex we will have a 1-deficient
or complete partite graph with at least 7 parts, all of which are
intrinsically linked.

\end{description}

\section{$2$-deficient graphs}

In this section we classify $2$-deficient graphs with respect to 
intrinsic linking and knotting. We also show that Adams' conjecture
holds for these graphs.

\Nt We will expand the notation that we created in the last section
by adding subscripts to the vertices. For example,
if we are removing two edges  between parts A and B of $K_{4,3,1}$,
there are 3 cases to be considered: deleting two edges that share a
vertex from part A ($K_{4,3,1}$-\{$(a_1,b_1),(a_1,b_2)$\}), deleting two
edges that share a vertex from part B
($K_{4,3,1}$-\{$(a_1,b_1),(a_2,b_1)$\}), and deleting two edges that
share no vertices ($K_{4,3,1}$-\{$(a_1,b_1),(a_2,b_2)$\}). Also, in some
cases, if we delete a particular edge we can then delete any other 
without affecting the classification of the resulting graph. For
example, if we delete edge (b,c) from $K_{4,3,1}$, we can then delete any
other and still have an intrinsically linked graph; we will denote graphs
obtained in this way by
$K_{4,3,1}$-\{$(b,c),e$\}. Furthermore, in some cases, we can delete any
2 edges and still have an intrinsically linked graph;
for example, $K_7$ will be
intrinsically linked no matter what 2 edges we delete. We denote these
graphs $K_7$-2e.

\subsection{Intrinsic linking}

\begin{table}[ht]

\begin{center}

\begin{tabular}{|l||c|c|c|c|} \hline

k          & 1    & 2      & 3                                & 4\\ \hline

linked     & 7-2e & 5,4-2e & 4,3,1-\{($b,c),e$\}                &
3,2,1,1-\{($b_1,c),(b_1,d$)\} \\ 

           &      &        & 4,3,1-\{($a_1,b_1),(a_1,b_2$)\}  &
3,2,1,1-\{$(b_1,c),(b_2,c)$\} \\ \cline{5-5}

           &      &        & 4,3,1-\{($a_1,b_1$),($a_1,c$)\}  &
4,2,1,1-\{$(b,c),e$\} \\ \cline{4-4}

           &      &        & 3,3,2-\{($a_1,b_1),(a_1,b_2$)\}  &
4,2,1,1-\{$(c,d),e$\} \\ 

           &      &        & 3,3,2-\{($a_1,c_1),(a_2,c_1$)\}  &
4,2,1,1-\{$(a_1,b_1),(a_1,b_2)$\} \\ 

           &      &        & 3,3,2-\{($a_1,c_1),(b_1,c_1$)\}  &
4,2,1,1-\{$(a_1,b_1),(a_1,c)$\} \\ 

           &      &        & 3,3,2-\{($a_1,b_1),(b_2,c_1$)\}  &
4,2,1,1-\{$(a_1,c),(a_1,d)$\} \\ \cline{5-5} 

           &      &        & 3,3,2-\{($a_1,c_1),(b_1,c_2$)\}  &
3,3,1,1-\{$(b,c),e$\} \\ 

           &      &        & 3,3,2-\{($a_1,c_1),(a_1,c_2$)\}   &
3,3,1,1-\{$(c,d),e$\} \\ \cline{4-4}

           &      &        & 4,2,2-\{$(b,c),e$\} &
3,3,1,1-\{$(a_1,b_1),(a_1,b_2)$\} \\ \cline{5-5}

           &      &        & 4,2,2-\{($a_1,b_1),(a_1,b_2$)\} &
3,2,2,1-\{$(b,c),e$\} \\ \cline{4-4}

           &      &        & 5,3,1-2e & 3,2,2,1-\{$(c,d),e$\} \\ 

           &      &        & 4,4,1-2e  & 3,2,2,1-\{$(a,d),e$\} \\ 

           &      &        & 4,3,2-2e   &
3,2,2,1-\{$(a_1,b_1),(a_1,b_2)$\} \\ 

           &      &        &  3,3,3-2e  &
3,2,2,1-\{$(a_1,b_1),(a_2,b_1)$\} \\ 

           &      &        &  5,2,2-2e &
3,2,2,1-\{$(a_1,b_1),(a_2,c_1)$\} \\ \cline{5-5}

           &      &        &   & 2,2,2,2-2e \\ 

           &      &        &                          & 5,2,1,1-2e \\ 

           &      &        &                          & 4,3,1,1-2e \\

           &      &        &                         & 4,2,2,1-2e \\

           &      &        &                         & 3,3,2,1-2e \\
\hline



not linked & 6-2e & 4,4-2e & 4,3,1-\{($a_1,b_1),(a_2,b_1$)\}  &
3,2,1,1-\{($a,b),e$\}\\ 

           &      & n,3-2e & 4,3,1-\{($a_1,b_1),(a_2,b_2$)\}  &
3,2,1,1-\{($a,c),e$\}\\ 

           &      &        & 4,3,1-\{($a_1,c),(a_2,c$)\}      &
3,2,1,1-\{($c,d),e$\}\\ 

           &      &        & 4,3,1-\{($a_1,b_1),(a_2,c$)\}    &
3,2,1,1-\{($b_1,c),(b_2,d)$\}\\ \cline{4-5} 

           &      &        & 3,3,2-\{($a_1,b_1),(a_2,b_2$)\}  &
4,2,1,1-\{$(a_1,b_1),(a_2,b_2)$\}\\ 

           &      &        & 3,3,2-\{($a_1,b_1),(b_1,c_1$)\}  &
4,2,1,1-\{$(a_1,c),(a_2,c)$\}\\ 

           &      &        & 3,3,2-\{($b_1,c_1),(b_2,c_2$)\}  &
4,2,1,1-\{$(a_1,b_1),(a_2,b_1$)\} \\ \cline{4-4}

           &      &        & 4,2,2-\{($a_1,b_1),(a_2,b_2$)\}  &
4,2,1,1-\{$(a_1,b_1),(a_2,c)$\}\\ 

           &      &        & 4,2,2-\{($a_1,b_1),(a_2,b_1$)\}  &
4,2,1,1-\{$(a_1,c),(a_2,d)$\}\\ \cline{5-5}

           &      &        & 4,2,2-\{($a_1,b_1),(a_2,c_1$)\}  &
3,3,1,1-\{$(a_1,b_1),(a_2,b_2$)\}\\ 

           &      &        & 4,2,2-\{($a_1,b_1),(a_1,c_1$)\}  &
3,2,2,1-\{$(a_1,b_1),(a_2,b_2$)\}\\ \cline{4-4}

           &      &        & 3,2,2-2e                         &
3,2,2,1-\{$(a_1,b_1),(a_1,c_1)$\}\\ \cline{5-5}

           &      &        & n,2,1-2e                         &
2,2,2,1-2e\\ 

           &      &        & 3,3,1-2e                         &
n,1,1,1-2e\\ \hline

\end{tabular}

\end{center}

\caption{Intrinsic Linking of 2 Deficient Graphs}\label{il2def1}

\end{table}

\begin{table}[ht]

\begin{center}

\begin{tabular}{|l||c|c|c|} \hline

k           & 5                                   &
6                                 & $\geq$7 \\ \hline

linked      & 2,2,1,1,1-\{$(b_1,c),(b_2,c)$\}     &
2,1,1,1,1,1-\{$(a_1,b),(a_1,c)$\}  & All \\ 

            & 2,2,1,1,1-\{$(a_1,d),(b_1,c)$\}     &
2,1,1,1,1,1-\{$(a_1,b),(a_2,b)$\}  &  \\ 

            & 2,2,1,1,1-\{$(b_1,c),(b_1,d)$\}     &
2,1,1,1,1,1-\{$(a_1,b),(c,d)$\} &  \\ \cline{2-2}

            & 3,1,1,1,1-\{$(b,c),(c,d)\}$         &
2,1,1,1,1,1-\{$(b,c),(c,d)$\} &  \\ \cline{2-3}

            & 4,1,1,1,1-\{$(b,c),e$\}             & 3,1,1,1,1,1-2e &  \\

            & 4,1,1,1,1-\{$(a_1,b),(a_1,c)$\}     & 2,2,1,1,1,1-2e &  \\
\cline{2-2}

            & 3,2,1,1,1-\{$(a,c),e$\}             &  &  \\

            & 3,2,1,1,1-\{$(b,c),e$\} &  &  \\

            & 3,2,1,1,1-\{$(c,d),e$\} &  &  \\

            & 3,2,1,1,1-\{$(a_1,b_1),(a_1,b_2)$\} &  &  \\

            & 3,2,1,1,1-\{$(a_1,b_1),(a_2,b_1)$\} &  &  \\ \cline{2-2}

            & 2,2,2,1,1-2e                        &  &  \\

            & 5,1,1,1,1-2e                        &  &  \\

            & 4,2,1,1,1-2e                        &  &  \\

            & 3,3,1,1,1-2e                        &  &  \\ \hline



not linked  & 2,2,1,1,1-\{$(a,b),e$\}             &
2,1,1,1,1,1-\{$(a_1,b),(a_2,c)$\}  & None \\ 

            & 2,2,1,1,1-\{$(c,d),e$\}             &
2,1,1,1,1,1-\{$(a_1,b),(b,c)$\} &  \\

            & 2,2,1,1,1-\{$(b_1,c),(b_2,d)$\}     &
2,1,1,1,1,1-\{$(b,c),(d,e)$\} &  \\ \cline{3-3}

            & 2,2,1,1,1-\{$(a_1,c),(b_1,c)$\}     & 1,1,1,1,1,1-2e &  \\
\cline{2-2}

            & 3,1,1,1,1-\{$(a,b),e$\}             &  &  \\

            & 3,1,1,1,1-\{$(b,c),(d,e)$\}         &  &  \\ \cline{2-2}

            & 4,1,1,1,1-\{$(a_1,b),(a_2,c)$\}     &  &  \\

            & 4,1,1,1,1-\{$(a_1,b),(a_2,b)$\}     &  &  \\  \cline{2-2}

            & 3,2,1,1,1-\{$(a_1,b_1),(a_2,b_2)$\} &  &  \\ \cline{2-2}

            & 2,1,1,1,1-2e                        &  &  \\ \hline


\end{tabular}

\end{center}

\caption{Intrinsic Linking of 2 Deficient Graphs (cont.)}\label{il2def2}

\end{table}

\begin{thm}

The 2-deficient graphs are classified with respect to intrinsic  linking
according to Tables \ref{il2def1} and \ref{il2def2}.
\end{thm}

\Pf

\noindent%
\underline{Linked:}



$K_7$-2e is intrinsically linked; there are 2 cases:
$K_7$-\{$(a,b),(b,c)$\} and $K_7$-\{$(a,b),(c,d)$\}.
For the first case,
we can simply delete vertex $b$ to get $K_6$; for the second case, we can
contract edge $(a,c)$ to get $K_6$. Notice that any 2-deficient graph
with 7 or more parts will have $K_7$-2e as a minor, and will therefore be
intrinsically linked.

$K_{5,4}$-2e is intrinsically linked; there are 3 cases. In each case, one
removed edge will be $(a_1,b_1)$. If we simply delete vertex $a_1$ we
will get $K_{4,4}$ or $K_{4,4}$-e.

By lemma \ref{combparts}, $K_{5,3,1}$-2e, $K_{4,4,1}$-2e, $K_{4,3,2}$-2e,
$K_{5,2,2}$-2e, $K_{5,2,1,1}$-2e, $K_{4,3,1,1}$-2e, $K_{4,2,2,1}$-2e,
$K_{3,3,2,1}$-2e, $K_{5,1,1,1,1}$-2e, $K_{4,2,1,1,1}$-2e, and
$K_{3,3,1,1,1}$-2e all have $K_{5,4}$-2e as a minor, and are therefore
intrinsically linked.

$K_{3,3,3}$-2e is intrinsically linked. Removing any 2 edges  will result
in a graph of the form $K_{3,3,3}$-\{$(a,b),e$\}. If we delete vertex
$a$, we will get $K_{3,3,2}$, or $K_{3,3,2}$-e, so $K_{3,3,3}$-2e is
intrinsically linked.

$K_{2,2,2,2}$-2e is intrinsically linked. All such graphs are of the form
$K_{2,2,2,2}$-\{$(a,b),e$\}. It them follows from lemma \ref{combparts}
that $K_{4,4}$-e is a minor of $K_{2,2,2,2}$-2e.

By lemma \ref{combparts}, $K_{2,2,2,1,1}$-2e and $K_{2,2,1,1,1,1}$-2e
have $K_{2,2,2,2}$-2e as a minor, and are therefore intrinsically linked.

$K_{3,1,1,1,1,1}$-2e is intrinsically linked. All such graphs are of the
form $K_{3,1,1,1,1,1}$-\{$(a,b),e$\} or $K_{3,1,1,1,1,1}$-\{$(b,c),e$\}.
The first case contains $K_{4,4}$-e as a minor by lemma \ref{combparts},
and the second case is equivalent to $K_{3,2,1,1,1}$-e, so
$K_{3,1,1,1,1,1}$-2e is intrinsically linked.

By lemma \ref{combparts}, $K_{4,3,1}$-\{$(b,c),e$\},
$K_{4,2,2}$-\{$(b,c),e$\}, $K_{4,2,1,1}$-\{$(b,c),e$\},
$K_{4,2,1,1}$-\{$(c,d),e$\}, $K_{3,3,1,1}$-\{$(b,c),e$\},
$K_{3,2,2,1}$-\{$(b,c),e$\}, $K_{3,2,2,1}$-\{$(a,d),e$\},
$K_{4,1,1,1,1}$-\{$(b,c),e$\}, $K_{3,2,1,1,1}$-\{$(a,c),e$\},
$K_{3,2,1,1,1}$-\{$(b,c),e$\}, and $K_{3,2,1,1,1}$-\{$(c,d),e$\} all have
$K_{4,4}$-e as a minor, and are therefore all intrinsically linked.

$K_{3,3,1,1}$-\{$(c,d),e$\}, and $K_{3,2,2,1}$-\{$(c,d),e$\} have
$K_{3,3,2}$-e as a minor by lemma \ref{combparts}, so they are
intrinsically linked.

$K_{4,3,1}$-\{$(a_1,b_1),(a_1,b_2)$\} and
$K_{4,3,1}$-\{$(a_1,b_1),(a_1,c)$\} both have $K_{3,3,1}$ as a minor;
simply delete vertex $a_1$. So both are intrinsically linked.

$K_{4,2,1,1}$-\{$(a_1,b_1),(a_1,b_2)$\},
$K_{4,2,1,1}$-\{$(a_1,b_1),(a_1,c)$\},
$K_{3,3,1,1}$-\{$(a_1,b_1),(a_1,b_2)$\},
$K_{3,2,2,1}$-\{$(a_1,b_1),(a_2,b_1)$\},
$K_{4,1,1,1,1}$-\{$(a_1,b),(a_1,c)$\},
$K_{3,2,1,1,1}$-\{$(a_1,b_1),(a_1,b_2)$\}, and
$K_{3,2,1,1,1}$-\{$(a_1,b_1),(a_2,b_1)$\} all contain
$K_{4,3,1}$-\{$(a_1,b_1),(a_1,b_2)$\} as a minor, so they are all
intrinsically linked.
(Notice that for the case of $K_{3,2,2,1}$ and
$K_{3,2,1,1,1}$-\{$(a_1,b_1),(a_2,b_1)$\}, once the parts are combined,
the labels ($a$ and $b$) of some of the vertices switch.)

$K_{3,3,2}$-\{$(a_1,b_1),(a_1,b_2)$\} and
$K_{3,3,2}$-\{$(a_1,b_1),(b_2,c_1)$\} have $K_{3,2,1,1}$ as a minor. To
see this, simply contract edge $(a_1,c_1)$. So both are intrinsically
linked.

By lemma \ref{combparts}, $K_{3,2,2,1}$-\{$(a_1,b_1),(a_1,b_2)$\} has
$K_{3,3,2}$-\{$(a_1,b_1),(a_1,b_2)$\} as a minor, so it is intrinsically
linked.

By lemma \ref{combparts}, $K_{3,2,2,1}$-\{$(a_1,b_1),(a_2,c_1)$\} has
$K_{3,3,2}$-\{$(a_1,b_1),(b_2,c_1)$\} as a minor. Notice that
$K_{3,3,2}$-\{$(a_1,b_1),(b_2,c_1)$\} and
$K_{3,3,2}$-\{$(b_1,a_1),(a_2,c_1)$\} are equivalent.

$K_{3,3,2}$-\{$(a_1,c_1),(a_2,c_1)$\} and
$K_{3,3,2}$-\{$(a_1,c_1),(b_1,c_1)$\} both have $K_{3,3,1}$ as a minor;
simply delete vertex $c_1$. So both are intrinsically linked.

$K_{3,3,2}$-\{$(a_1,c_1),(b_1,c_2)$\} is intrinsically linked. Delete the
remaining edges between $c_1$ and A and between $c_2$ and B. The result
is $K_{4,4}-e$ where the missing edge is between vertices $c_1$ and $c_2$.

\begin{figure}[ht]
\begin{center}
\includegraphics[scale=0.3]{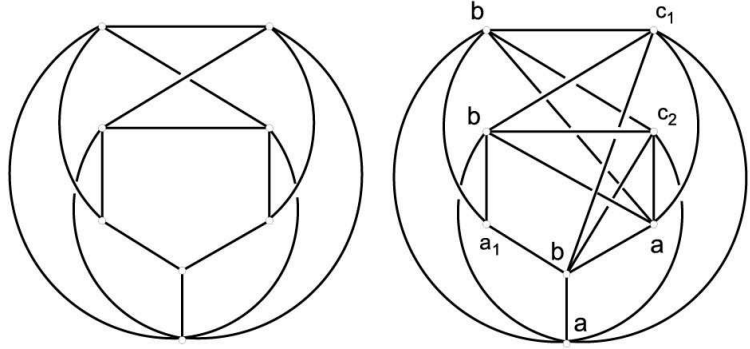}

\caption{$K_{3,3,2}-\{(a_1,c_1),(a_1,c_2)\}$ is Intrinsically
Linked}\label{k332a1c1a1c2}
\end{center}
\end{figure}

$K_{3,3,2}-\{(a_1,c_1),(a_1,c_2)\}$ is intrinsically linked. In
figure \ref{k332a1c1a1c2}, at left is a Petersen graph, and at
right is that same graph with four edges added. As we can see, this
is $K_{3,3,2}-\{(a_1,c_1),(a_1,c_2)\}$.

$K_{4,2,2}$-\{$(a_1,b_1),(a_1,b_2)$\} is intrinsically linked; contract
edge $(a_1,c_1)$ to get $K_{3,2,1,1}$

$K_{4,2,1,1}$-\{$(a_1,c),(a_1,d)$\} has
$K_{4,2,2}$-\{$(a_1,b_1),(a_1,b_2)$\} as a minor by lemma \ref{combparts}.

$K_{3,2,1,1}$-\{($b_1,c),(b_2,c)$\} is equivalent to $K_{3,3,1}$ which is
intrinsically linked.

$K_{3,2,1,1}$-\{$(b_1,c),(b_1,d)$\} is a Petersen graph; it is obtained by
a triangle-Y exchange of $K_6$.

$K_{2,2,1,1,1}$-\{$(b_1,c),(b_2,c)$\} is equivalent to $K_{3,2,1,1}$, so
it is intrinsically linked.

$K_{2,2,1,1,1}$-\{$(a_1,d),(b_1,c)$\} has $K_{3,3,1}$ as a minor by lemma
\ref{combparts}.

$K_{2,2,1,1,1}$-\{$(b_1,c),(b_1,d)$\} has $K_6$ as a minor; simply
contract edge $(a_1,b_1)$.

$K_{3,1,1,1,1}$-\{$(b,c),(c,d)$\} has $K_{3,3,1}$ as a minor by lemma
\ref{combparts}.

$K_{2,1,1,1,1,1}$-\{$(a_1,b),(a_1,c)$\} has $K_6$ as a minor; simply
delete vertex $a_1$.

$K_{2,1,1,1,1,1}$-\{$(a_1,b),(a_2,b)$\} is equivalent to $K_{3,1,1,1,1}$.

$K_{2,1,1,1,1,1}$-\{$(a_1,b),(c,d)$\} is equivalent to
$K_{2,2,1,1,1}$-$(a,c)$

$K_{2,1,1,1,1,1}$-\{$(b,c),(c,d)$\} has $K_{3,2,1,1}$ as a minor.

\medskip

\noindent%
\underline{Not Linked:}

The following graphs are already known not to be intrinsically linked if
we add one more edge, so they are naturally not intrinsically linked:
$K_6$-2e (or $K_{1,1,1,1,1,1}$-2e), $K_{n,3}$-2e, $K_{3,2,2}$-2e,
$K_{n,2,1}$-2e, $K_{3,3,1}$-2e, $K_{2,2,2,1}$-2e, $K_{n,1,1,1}$-2e,
$K_{3,2,1,1}$-\{$(a,b),e$\}, $K_{3,2,1,1}$-\{$(a,c),e$\},
$K_{3,2,1,1}$-\{$(c,d),e$\}, $K_{2,2,1,1,1}$-\{$(a,b),e$\},
$K_{2,2,1,1,1}$-\{$(c,d),e$\}, $K_{3,1,1,1,1}$-\{$(a,b),e$\},
$K_{2,1,1,1,1}$-2e.

$K_{2,1,1,1,1,1}$-\{$(b,c),(d,e)$\},
$K_{2,1,1,1,1,1}$-\{$(a_1,b),(b,c)$\},
$K_{2,1,1,1,1,1}$-\{$(a_1,b),(a_2,c)$\},
$K_{4,1,1,1,1}$-\{$(a_1,b),(a_2,b)$\},
$K_{4,1,1,1,1}$-\{$(a_1,b),(a_2,c)$\},
$K_{3,2,1,1,1}$-\{$(a_1,b_1),(a_2,b_2)$\},
$K_{2,2,1,1,1}$-\{$(a_1,c),(b_1,c)$\},
$K_{2,2,1,1,1}$-\{$(b_1,c),(b_2,d)$\},
$K_{3,2,2,1}$-\{$(a_1,b_1),(a_1,c_1)$\},
$K_{3,2,2,1}$-\{$(a_1,b_1),(a_2,b_2)$\},
$K_{3,3,1,1}$-\{$(a_1,b_1),(a_2,b_2)$\},
$K_{4,2,1,1}$-\{$(a_1,b_1),(a_2,b_2)$\},
$K_{4,2,1,1}$-\{$(a_1,c),(a_2,c)$\},
$K_{4,2,1,1}$-\{$(a_1,b_1),(a_2,b_1)$\},
$K_{4,2,1,1}$-\{$(a_1,b_1),(a_2,c)$\},
$K_{4,3,1}$-\{$(a_1,b_1),(a_2,b_1)$\}, and
$K_{4,3,1}$-\{$(a_1,b_1),(a_2,b_2)$\} all have 1 vertex connected to
every other. In each of these, if we delete that vertex we get a planar
graph; so, by lemma \ref{sachs}, none are intrinsically linked.

$K_{3,1,1,1,1}$-\{$(b,c),(d,e)$\} is equivalent to $K_{3,2,2}$, so it is
not intrinsically linked.

$K_{4,2,1,1}$-\{$(a_1,c),(a_2,d)$\} is a minor of
$K_{4,1,1,1,1}$-\{$(a_1,b),(a_2,c)$\} by lemma \ref{combparts}; simply
combine parts D and E.

$K_{3,2,1,1}$-\{$(b_1,c),(b_2,d)$\} is a minor of
$K_{3,1,1,1,1}$-\{$(b,c),(d,e)$\} by lemma \ref{combparts}; simply
combine parts C and D.

$K_{4,2,2}$-\{$(a_1,b_1),(a_2,b_2)$\},
$K_{4,2,2}$-\{$(a_1,b_1),(a_2,b_1)$\}, and
$K_{4,2,2}$-\{$(a_1,b_1),(a_2,c_1)$\}, are minors of
$K_{4,2,1,1}$-\{$(a_1,b_1),(a_2,b_2)$\},
$K_{4,2,1,1}$-\{$(a_1,b_1),(a_2,b_1)$\}, and
$K_{4,2,1,1}$-\{$(a_1,b_1),(a_2,c)$\} respectively by lemma
\ref{combparts}.

$K_{4,2,2}-\{(a_1,b_1),(a_1,c_1)\}$ is a minor of 
$K_{3,2,2,1}-\{(a_1,b_1),(a_1,c_1)\}$ and is, therefore, not
intrinsically linked.

$K_{3,3,2}$-\{$(a_1,b_1),(a_2,b_2)$\} is a minor of
$K_{3,3,1,1}$-\{$(a_1,b_1),(a_2,b_2)$\}.

$K_{3,3,2}$-\{$(b_1,c_1),(b_2,c_2)$\} is a minor of
$K_{3,2,2,1}$-\{$(a_1,b_1),(a_2,b_2)$\}.

$K_{3,3,2}$-\{$(a_1,b_1),(b_1,c_1)$\} is not intrinsically linked by the
corollary to lemma~\ref{sachs}; if you delete vertex $a_2$, the result is
a planar graph.

$K_{4,3,1}$-\{$(a_1,c),(a_2,c)$\} and $K_{4,3,1}$-\{$(a_1,b_1),(a_2,c)$\}
are minors of $K_{4,2,1,1}$-\{$(a_1,c),(a_2,c)$\} and
$K_{4,2,1,1}$-\{$(a_1,b_1),(a_2,c)$\} respectively by lemma
\ref{combparts}, and therefore are not intrinsically linked.

$K_{4,4}$-2e is a minor of $K_{4,4}$-e, which is minor minimally
intrinsically linked, so it is not intrinsically linked. \qed

\subsection{Intrinsic knotting}

\footnotesize

\begin{table}[ht]

\begin{center}

\begin{tabular}{|l||c|c|c|c|c|} \hline

k           & 1    & 2      & 3                                 & 4&5\\
\hline

knotted     & 8-2e & 5,5-2e & 3,3,3-\{($a_1,b_1$),($a_1,b_2$)\} &
3,2,2,2-\{$(b,c),e$\}& 

3,2,1,1,1-\{($b_1,c$),($b_1,d$)\} \\

            &      &        & 3,3,3-\{($a_1,b_1$),($b_2,c_1$)\} &
3,2,2,2-\{$(a_1,b_1),(a_1,b_2)$\} & 

3,2,1,1,1-\{($b_1,c$),($b_2,c$)\} \\ \cline{4-4} \cline{6-6}

            &      &        & 4,4,1-\{($a_1,c$),($b_1,c$)\}     &
3,2,2,2-\{$(a_1,b_1),(a_2,b_1)$\} & 

4,2,1,1,1-\{$(c,d),e$\}  \\

            &      &        & 4,4,1-\{($a_1,b_1$),($b_1,c$)\}   &
3,2,2,2-\{$(a_1,b_1),(a_2,c_1)$\} & 

4,2,1,1,1-\{$(b,c),e$\}  \\ \cline{4-5}

            &      &        & 4,3,2-\{($a_1,b_1$),($a_1,b_2$)\} &
4,2,2,1-\{$(b,c),e$\} &

4,2,1,1,1-\{($a_1,b_1$),($a_1,c$)\} \\

            &      &        & 4,3,2-\{($a_1,c_1$),($b_1,c_2$)\} &
4,2,2,1-\{$(c,d),e$\} &

4,2,1,1,1-\{($a_1,b_1$),($a_1,b_2$)\} \\

            &      &        & 4,3,2-\{($a_1,c_1$),($b_1,c_1$)\} &
4,2,2,1-\{($a_1,b_1$),($a_1,b_2$)\} &

4,2,1,1,1-\{($a_1,c$),($a_1,d$)\} \\ \cline{6-6}

            &      &        & 4,3,2-\{($b_1,c_1$),($b_2,c_1$)\} &
4,2,2,1-\{($a_1,b_1$),($a_1,d$)\} &

3,3,1,1,1-\{$(b,c),e$\}\\  \cline{5-5}

            &      &        & 4,3,2-\{($a_1,c_1$),($a_1,c_2$)\} &
3,3,2,1-\{$(c,d),e$\} &

3,3,1,1,1-\{$(c,d),e$\}\\

            &      &        & 4,3,2-\{($a_1,b_1$),($b_2,c_1$)\} &
3,3,2,1-\{$(a,d),e$\} &

3,3,1,1,1-\{$(a_1,b_1),(a_2,b_1)$\}\\ \cline{6-6}

            &      &        & 4,3,2-\{($b_1,c_1$),($b_1,c_2$)\} &
3,3,2,1-\{($a_1,b_1$),($a_1,b_2$)\} &

3,2,2,1,1-\{$(a,d),e$\}\\ \cline{4-4}

            &      &        & 5,4,1-2e                          &
3,3,2,1-\{($b_1,c_1$),($b_1,c_2$)\} &

3,2,2,1,1-\{$(c,d),e$\}\\

            &      &        & 5,3,2-2e                          &
3,3,2,1-\{($a_1,c_1$),($b_1,c_1$)\} &

3,2,2,1,1-\{$(b,c),e$\}\\

            &      &        & 4,3,3-2e                          &
3,3,2,1-\{($a_1,b_1$),($a_2,c_1$)\} &

3,2,2,1,1-\{$(d,e),e$\}\\

            &      &        & 4,4,2-2e              &
3,3,2,1-\{($a_1,c_1$),($a_2,c_1$)\}
&3,2,2,1,1-\{$(a_1,b_1),(a_1,b_2)$\}\\ 

            &      &        &           &
3,3,2,1-\{($a_1,c_1$),($b_1,c_2$)\}
&3,2,2,1,1-\{$(a_1,b_1),(a_2,b_1)$\}\\  \cline{5-5}

            &      &        &         & 4,3,1,1-\{$(b,c),e$\}
&3,2,2,1,1-\{$(a_1,b_1),(a_2,c_1)$\}\\ \cline{6-6}

            &      &        &                                   &
4,3,1,1-\{$(c,d),e$\} &2,2,2,2,1-2e\\ 

            &      &        &                                   &
4,3,1,1-\{($a_1,b_1$),($a_1,b_2$)\} &5,2,1,1,1-2e\\ 

            &      &        &                                   &
4,3,1,1-\{($a_1,b_1$),($a_1,c$)\} &4,3,1,1,1-2e\\ 

            &      &        &                                   &
4,3,1,1-\{($a_1,c$),($a_1,d$)\} &4,2,2,1,1-2e\\ \cline{5-5}

            &      &        &                                   &
4,2,2,2-2e &3,3,2,1,1-2e \\ 


            &      &        &                                   &
3,3,2,2-2e &\\  

            &      &        &                                   &
5,2,2,1-2e & \\ 

            &      &        &                                   &
4,3,2,1-2e & \\ 

            &      &        &                                   &
3,3,3,1-2e & \\ 

            &      &        &                                   &
5,3,1,1-2e & \\ 

            &      &        &                                   &
4,4,1,1-2e & \\ \hline


not knotted & 7-2e & n,4-2e & 3,3,3-\{($a_1,b_1$),($b_1,c_1$)\} &
3,2,2,2-\{($a_1,b_1$),($a_1,c_1$)\} & 3,2,1,1,1-\{$(a,b),e$\}\\

            &      &        & 3,3,3-\{($a_1,b_1$),($a_2,b_2$)\} &
3,2,2,2-\{($a_1,b_1$),($a_2,b_2$)\} &3,2,1,1,1-\{$(a,c),e$\}\\ \cline{4-5}

            &      &        & 4,4,1-\{($a_1,b_1$),($a_1,b_2$)\} &
4,2,2,1-\{($a_1,b_1$),($a_2,b_1$)\} &3,2,1,1,1-\{$(c,d),e$\}\\

            &      &        & 4,4,1-\{($a_1,b_1$),($a_2,b_2$)\} &
4,2,2,1-\{($a_1,d$),($a_2,d$)\} &

3,2,1,1,1-\{$(b_1,c),(b_2,d)$\}\\ \cline{6-6}

            &      &        & 4,4,1-\{($a_1,c$),($a_2,c$)\}     &
4,2,2,1-\{($a_1,b_1$),($a_2,b_2$)\}
&4,2,1,1,1-\{($a_1,b_1$),($a_2,b_2$)\}\\

            &      &        & 4,4,1-\{($a_1,b_1$),($a_2,c$)\}   &
4,2,2,1-\{($a_1,b_1$),($a_2,c_1$)\}&4,2,1,1,1-\{($a_1,b_1$),($a_2,c$)\}\\
\cline{4-4}

            &      &        & 4,3,2-\{($a_1,b_1$),($a_2,b_1$)\} &
4,2,2,1-\{($a_1,b_1$),($a_2,d$)\} &4,2,1,1,1-\{($a_1,b_1$),($a_2,b_1$)\}\\

            &      &        & 4,3,2-\{($a_1,b_1$),($a_2,b_2$)\} &
4,2,2,1-\{($a_1,b_1$),($a_1,c_1$)\} &4,2,1,1,1-\{($a_1,c$),($a_2,c$)\}\\
\cline{5-5}

            &      &        & 4,3,2-\{($a_1,c_1$),($a_2,c_1$)\} &
3,3,2,1-\{($a_1,b_1$),($b_1,c_1$)\} &4,2,1,1,1-\{($a_1,c$),($a_2,d$)\}\\
\cline{6-6}

            &      &        & 4,3,2-\{($a_1,b_1$),($b_1,c_1$)\} &
3,3,2,1-\{($a_1,b_1$),($a_2,b_2$)\}
&3,3,1,1,1-\{($a_1,b_1$),($a_2,b_2$)\}\\ \cline{6-6}

            &      &        & 4,3,2-\{($a_1,b_1$),($a_2,c_1$)\} &
3,3,2,1-\{($b_1,c_1$),($b_2,c_2$)\} &
3,2,2,1,1-\{($a_1,b_1$),($a_1,c_1$)\}\\ \cline{5-5}

            &      &        & 4,3,2-\{($b_1,c_1$),($b_2,c_2$)\} &
4,3,1,1-\{($a_1,b_1$),($a_2,b_1$)\}
&3,2,2,1,1-\{($a_1,b_1$),($a_2,b_2$)\} \\ \cline{6-6}

            &      &        & 4,3,2-\{($a_1,c_1$),($a_2,c_2$)\} &
4,3,1,1-\{($a_1,b_1$),($a_2,b_2$)\} &2,2,2,1,1-2e \\

            &      &        & 4,3,2-\{($a_1,b_1$),($a_1,c_1$)\} &
4,3,1,1-\{($a_1,b_1$),($a_2,c$)\} &n,1,1,1,1-2e\\ \cline{4-4}

            &      &        & 3,3,2-2e                          &
4,3,1,1-\{($a_1,c$),($a_2,d$)\} &\\

            &      &        & n,2,2-2e                          &
4,3,1,1-\{($a_1,c$),($a_2,c$)\} &\\ \cline{5-5}

            &      &        & n,3,1-2e                          &
3,3,1,1-2e &\\

            &      &        &                                   &
2,2,2,2-2e & \\ 

            &      &        &                                   &
3,2,2,1-2e &\\ 

            &      &        &                                   &
n,2,1,1-2e &\\  \hline


\end{tabular}

\end{center}

\caption{Intrinsic Knotting of 2 Deficient Graphs}\label{ik2def1}

\end{table}

\normalsize

\begin{table}[ht]

\begin{center}

\begin{tabular}{|l||c|c|c|} \hline

k             & 6                                   & 7 & $\geq$8 \\
\hline

knotted       & 2,2,1,1,1,1-\{($b_1,c$),($b_2,c$)\} &
2,1,1,1,1,1,1-\{($a_1,b$),($a_1,c$)\} & all \\ 

              & 2,2,1,1,1,1-\{($b_1,c$),($b_1,d$)\} &
2,1,1,1,1,1,1-\{($a_1,b$),($a_2,b$)\} &  \\ 

              & 2,2,1,1,1,1-\{($a_1,d$),($b_1,c$)\}     &
2,1,1,1,1,1,1-\{($a_1,b$),($c,d$)\} &  \\ \cline{2-2}

              & 3,1,1,1,1,1-\{($b,c$),($c,d$)\}           &
2,1,1,1,1,1,1-\{($b,c$),($c,d$)\} &  \\ \cline{2-3}

              & 4,1,1,1,1,1-\{$(b,c),e$\} & 3,1,1,1,1,1,1-2e &  \\ 

              & 4,1,1,1,1,1-\{($a_1,b$),($a_1,c$)\}      &
2,2,1,1,1,1,1-2e &  \\ \cline{2-2}

              & 3,2,1,1,1,1-\{$(a,c),e$\}   &  &  \\ 

              &  3,2,1,1,1,1-\{$(b,c),e$\}  &  &  \\ 

              & 3,2,1,1,1,1-\{$(c,d),e$\}    &  &  \\ 

              & 3,2,1,1,1,1-\{$(a_1,b_1),(a_1,b_2)$\}    &  &  \\ 

              & 3,2,1,1,1,1-\{$(a_1,b_1),(a_2,b_1)$\} &  &  \\ \cline{2-2}

              & 2,2,2,1,1,1-2e &  &  \\ 

              & 5,1,1,1,1,1-2e &  &  \\ 

              & 4,2,1,1,1,1-2e &  &  \\ 

              & 3,3,1,1,1,1-2e &  &  \\ \hline

not knotted   & 2,2,1,1,1,1-\{$(a,b),e$\}     &
2,1,1,1,1,1,1-\{($a_1,b$),($b,c$)\} & none \\ 

              &  2,2,1,1,1,1-\{$(c,d),e$\}     &
2,1,1,1,1,1,1-\{($b,c$),($d,e$)\} &  \\ 

              & 2,2,1,1,1,1-\{$(b_1,c),(b_2,d)$\} &
2,1,1,1,1,1,1-\{($a_1,b$),($a_2,c$)\} &  \\ \cline{3-3}

              & 2,2,1,1,1,1-\{$(a_1,c),(b_1,c)$\} & 1,1,1,1,1,1,1-2e & 
\\  \cline{2-2}

              & 3,1,1,1,1,1-\{$(a,b),e$\} &  &  \\ 

              & 3,1,1,1,1,1-\{$(b,c),(d,e)$\} &  &  \\ \cline{2-2}

              & 4,1,1,1,1,1-\{($a_1,b$),($a_2,b$)\} &  &  \\ 

              & 4,1,1,1,1,1-\{($a_1,b$),($a_2,c$)\} &  &  \\ \cline{2-2}

              & 3,2,1,1,1,1-\{($a_1,b_1$),($a_2,b_2$)\} &  &  \\ \cline{2-2}

              & 2,1,1,1,1,1-2e &  &  \\  \hline

\end{tabular}

\end{center}

\caption{Intrinsic Knotting of 2 Deficient Graphs (cont.)}\label{ik2def2}

\end{table}

\begin{thm}

The 2-deficient graphs are classified with respect to intrinsic  knotting
according to Tables \ref{ik2def1} and \ref{ik2def2}.
\end{thm}

\Pf

\noindent%
\underline{Knotted:}

$K_8$-2e has 2 cases, $K_8$-\{$(a,b),(b,c)$\} and
$K_8$-\{$(a,b),(c,d)$\}. In the first case, delete vertex $b$ to get
$K_7$; in the second, contract edge $(a,c)$ to get $K_7$. So $K_8$-2e is
intrinsically knotted. Note that any 2-deficient graph with 8 or more
parts will have $K_8$-2e as a minor.

\begin{figure}[ht]
\begin{center}
\includegraphics[scale=0.35]{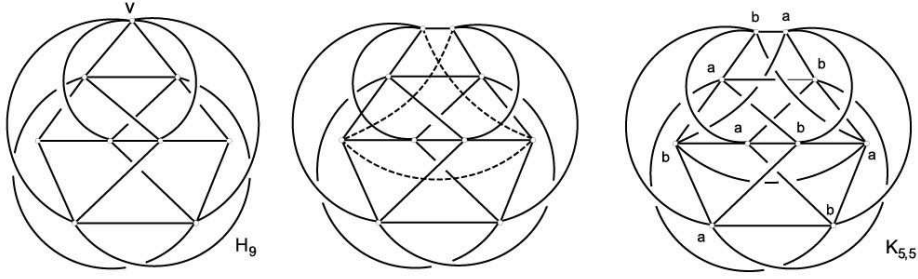}

\caption{$H_9$ is a Minor of $K_{5,5}$-2e}\label{h9tok55}
\end{center}
\end{figure}

$K_{5,5}$-2e has 2 cases. As shown in figure \ref{h9tok55}, by splitting
vertex $v$ and then adding the dashed lines, it can be seen that, in either case,
$K_{5,5}$-2e has $H_9$ (from \cite{KS}) as a minor.

By lemma \ref{combparts}, $K_{5,4,1}$-2e, $K_{5,3,2}$-2e,
$K_{3,3,2,2}$-2e, $K_{5,2,2,1}$-2e, $K_{4,3,2,1}$-2e, $K_{5,3,1,1}$-2e,
$K_{4,4,1,1}$-2e, $K_{5,2,1,1,1}$-2e, $K_{4,3,1,1,1}$-2e,
$K_{4,2,2,1,1}$-2e, $K_{3,3,2,1,1}$-2e, $K_{5,1,1,1,1,1}$-2e,
$K_{4,2,1,1,1,1}$-2e, and $K_{3,3,1,1,1,1}$-2e all contain $K_{5,5}$-2e
as a minor and are therefore intrinsically knotted.

$K_{4,3,3}$-2e is intrinsically knotted. In all cases, there is one
removed edge connected to vertex $b_1$; delete that vertex to get
$K_{4,3,2}$-e.

$K_{3,3,3,1}$-2e has $K_{4,3,3}$-2e as a minor by lemma \ref{combparts},
so it is intrinsically knotted.

$K_{4,4,2}$-2e is intrinsically knotted. In all cases, there is one removed
edge connected to vertex $a_1$; delete vertex $a_1$ to get
$K_{4,3,2}$-e.

By lemma \ref{combparts}, $K_{4,2,2,2}$-2e has $K_{4,4,2}$-2e as a minor.

$K_{2,2,2,2,1}$-2e is intrinsically knotted. In all cases, there will be
an edge missing between A and B or between D and E. If we combine parts A
and B in the first case, we get $K_{4,2,2,1}$-e or $K_{4,2,2,1}$, and if we combine
parts D and E in the second case, we get $K_{3,2,2,2}$-e or $K_{3,2,2,2}$. So, in
either case,
$K_{2,2,2,2,1}$-2e has an intrinsically knotted minor.

By lemma \ref{combparts}, $K_{2,2,2,1,1,1}$-2e and $K_{2,2,1,1,1,1,1}$-2e
have $K_{2,2,2,2,1}$-2e as a minor.

$K_{3,1,1,1,1,1,1}$-2e is intrinsically knotted. There is either an edge
missing between parts A and B or between parts B and C. In the first
case, we can combine A and B to get $K_{4,1,1,1,1,1}$-e or $K_{4,1,1,1,1,1}$, and
in the second, we can combine parts B and C to get $K_{3,2,1,1,1,1}$-e
or $K_{3,2,1,1,1}$.

By lemma \ref{combparts}, $K_{3,2,2,2}$-\{$(b,c),e$\},
$K_{4,2,2,1}$-\{$(c,d),e$\}, $K_{3,3,2,1}$-\{$(a,d),e$\},
$K_{4,3,1,1}$-\{$(c,d),e$\}, $K_{4,2,1,1,1}$-\{$(c,d),e$\},
$K_{4,2,1,1,1}$-\{$(b,c),e$\}, $K_{3,3,1,1,1}$-\{$(b,c),e$\},
$K_{3,3,1,1,1}$-\{$(c,d),e$\}, $K_{3,2,2,1,1}$-\{$(a,d),e$\},
$K_{3,2,2,1,1}$-\{$(c,d),e$\}, $K_{3,2,2,1,1}$-\{$(b,c),e$\},
$K_{3,2,2,1,1}$-\{$(d,e),e$\},
$K_{4,1,1,1,1,1}$-\{$(b,c),e$\}, $K_{3,2,1,1,1,1}$-\{$(a,c),e$\},
$K_{3,2,1,1,1,1}$-\{$(b,c),e$\}, $K_{3,2,1,1,1,1}$-\{$(c,d),e$\}, all
have $K_{4,3,2}$-e as a minor, and are therefore all intrinsically
knotted.

By lemma \ref{combparts}, $K_{4,2,2,1}$-\{$(b,c),e$\} and
$K_{4,3,1,1}$-\{$(b,c),e$\} have $K_{4,4,1}$-e as a minor.

$K_{3,3,2,1}$-\{$(c,d),e$\} has $K_{3,3,3}-e$ as a minor by lemma
\ref{combparts}.

As shown in figure \ref{h9tok333}, $H_9$ is a minor of $K_{3,3,3}$. If we
simply add fewer edges, we see that
$K_{3,3,3}$-\{$(a_1,b_1),(a_1,b_2)$\} and
$K_{3,3,3}$-\{$(a_1,b_1),(b_2,c_1)$\} are intrinsically knotted.





$K_{3,3,2,1}$-\{$(a_1,b_1),(a_1,b_2)$\},
$K_{3,3,2,1}$-\{$(b_1,c_1),(b_1,c_2)$\},
$K_{3,3,1,1,1}$-\{$(a_1,b_1),(a_2,b_1)$\},
$K_{3,2,2,1,1}$-\{$(a_1,b_1),(a_1,b_2)$\},
$K_{3,2,2,1,1}$-\{$(a_1,b_1),(a_2,b_1)$\},
$K_{3,2,1,1,1,1}$-\{$(a_1,b_1),(a_1,b_2)$\}, 
$K_{3,2,1,1,1,1}$-\{$(a_1,b_1),(a_2,b_1)$\}, and
$K_{3,3,2,1}$-\{$(a_1,c_1),(a_2,c_1)$\}  have
$K_{3,3,3}$-\{$(a_1,b_1),(a_1,b_2)$\} as a minor by lemma \ref{combparts}.

$K_{3,2,2,1,1}$-\{$(a_1,b_1),(a_2,c_1)$\},
$K_{3,3,2,1}$-\{$(a_1,b_1),(a_2,c_1)$\}, and
$K_{3,3,2,1}$-\{$(a_1,c_1),(b_1,c_2)$\} have
$K_{3,3,3}$-\{$(a_1,b_1),(b_2,c_1)$\} as a minor by lemma \ref{combparts}.

As shown in figure \ref{h9tok441}, $H_9$ is a minor of $K_{4,4,1}$. If we
simply add fewer edges, we can see that $K_{4,4,1}$-\{$(a_1,c),(b_1,c)$\}
and $K_{4,4,1}$-\{$(a_1,b_1),(b_1,c)$\} are intrinsically knotted.

$K_{4,2,2,1}$-\{$(a_1,b_1),(a_1,d)$\},
$K_{4,3,1,1}$-\{$(a_1,b_1),(a_1,c)$\},
$K_{4,2,1,1,1}$-\{$(a_1,b_1),(a_1,c)$\},
$K_{4,2,1,1,1}$-\{$(a_1,c),(a_1,d)$\},
$K_{4,1,1,1,1,1}$-\{$(a_1,b),(a_1,c)$\}, and
$K_{4,3,1,1}$-\{$(a_1,c),(a_1,d)$\} have
$K_{4,4,1}$-\{$(a_1,b_1),(b_1,c)$\} as a minor by lemma \ref{combparts}.

$H_9$ is a minor of $K_{4,3,2}$ as shown in figure \ref{h9tok432}. If we
add 2 fewer edges, we can see that $K_{4,3,2}$-\{$(a_1,c_1),(b_1,c_2)$\},
$K_{4,3,2}$-\{$(a_1,c_1),(b_1,c_1)$\},
$K_{4,3,2}$-\{$(b_1,c_1),(b_2,c_1)$\}, and
$K_{4,3,2}$-\{$(a_1,c_1),(a_1,c_2)$\} have $H_9$ as a minor, and are
therefore intrinsically knotted.

As shown in figure \ref{k3311tok432}, $K_{3,3,1,1}$ is a minor of
$K_{4,3,2}$. If we add 2 fewer edges, we can see that
$K_{4,3,2}$-\{$(a_1,b_1),(a_1,b_2)$\} has $K_{3,3,1,1}$ as a minor, so it
is intrinsically knotted.

$K_{4,3,2}$-\{$(a_1,b_1),(b_2,c_1)$\} is intrinsically knotted; if we
contract edge $(a_1,c_1)$ we get $K_{3,3,1,1}$.

$K_{4,3,2}$-\{$(b_1,c_1),(b_1,c_2)$\} is intrinsically knotted; it has
$B_9$ from table \ref{k3311} as a minor.

$K_{3,2,2,2}$-\{$(a_1,b_1),(a_1,b_2)$\} has
$K_{4,3,2}$-\{$(b_1,c_1),(b_1,c_2)$\} as a minor by lemma \ref{combparts}.

$K_{3,2,2,2}$-\{$(a_1,b_1),(a_2,b_1)$\},
$K_{4,2,2,1}$-\{$(a_1,b_1),(a_1,b_2)$\},
$K_{4,3,1,1}$-\{$(a_1,b_1),(a_1,b_2)$\}, and
$K_{4,2,1,1,1}$-\{$(a_1,b_1),(a_1,b_2)$\}  have
$K_{4,3,2}$-\{$(a_1,b_1),(a_1,b_2)$\} as a minor by lemma \ref{combparts}.

$K_{3,3,2,1}$-\{$(a_1,c_1),(b_1,c_1)$\} has
$K_{4,3,2}$-\{$(a_1,c_1),(b_1,c_1)$\} as a minor by lemma \ref{combparts}.

$K_{3,2,2,2}$-\{$(a_1,b_1),(a_2,c_1)$\} has
$K_{4,3,2}$-\{$(a_1,b_1),(b_2,c_1)$\} as a minor by lemma \ref{combparts}.

$K_{3,2,1,1,1}$-\{$(b_1,c),(b_2,c)$\} is equivalent to $K_{3,3,1,1}$.

$K_{3,2,1,1,1}$-\{$(b_1,c),(b_1,d)$\} has $H_8$ from \cite{KS}
as a minor.

$K_{2,2,1,1,1,1}$-\{$(b_1,c),(b_1,d)$\},
$K_{3,1,1,1,1,1}$-\{$(b,c),(c,d)$\}, and
$K_{2,1,1,1,1,1,1}$-\{$(a_1,b),(a_1,c)$\} have
$K_{3,2,1,1,1}$-\{$(b_1,c),(b_1,d)$\} as a minor by lemma \ref{combparts}.

$K_{2,2,1,1,1,1}$-\{$(b_1,c),(b_2,c)$\},
$K_{2,1,1,1,1,1,1}$-\{$(a_1,b),(a_2,b)$\}, and
$K_{2,1,1,1,1,1,1}$-\{$(b,c),(c,d)$\} have
$K_{3,2,1,1,1}$-\{$(b_1,c),(b_2,c)$\} as a minor by lemma \ref{combparts}.

$K_{2,2,1,1,1,1}$-\{$(a_1,d),(b_1,c)$\} has $K_{3,3,1,1}$ as a minor.

$K_{2,1,1,1,1,1,1}$-\{$(a_1,b),(c,d)$\} has 2 vertices connected to every
other, if you delete those, you get a non-planar graph, so by lemma
\ref{flemming}, it is intrinsically knotted.

\medskip

\noindent%
\underline{Not Knotted:}

Any graph that was not knotted with 1 edge removed will clearly not be
knotted with 2 edges removed. The following graphs fit that description:
$K_7$-2e ($= K_{1,1,1,1,1,1,1}$-2e), $K_{n,4}$-2e, $K_{3,3,2}$-2e,
$K_{n,2,2}$-2e, $K_{n,3,1}$-2e, $K_{3,3,1,1}$-2e, $K_{2,2,2,2}$-2e,
$K_{3,2,2,1}$-2e, $K_{n,2,1,1}$-2e, $K_{3,2,1,1,1}$-\{$(a,b),e$\},
$K_{3,2,1,1,1}$-\{$(a,c),e$\}, $K_{3,2,1,1,1}$-\{$(c,d),e$\},
$K_{2,2,2,1,1}$-2e, $K_{n,1,1,1,1}$-2e, $K_{2,2,1,1,1,1}$-\{$(a,b),e$\}, 
$K_{2,2,1,1,1,1}$-\{$(c,d),e$\}, $K_{3,1,1,1,1,1}$-\{$(a,b),e$\}, and
$K_{2,1,1,1,1,1}$-2e.

If a graph has 2 vertices connected to every other vertex and to one
another and the deletion of those vertices results in a planar graph,
lemma \ref{flemming} states that the original graph is not intrinsically
knotted. The following graphs are of that form:
$K_{4,3,1,1}$-\{$(a_1,b_1),(a_2,b_1)$\},
$K_{4,3,1,1}$-\{$(a_1,b_1),(a_2,b_2)$\},
$K_{4,2,1,1,1}$-\{$(a_1,b_1),(a_2,b_2)$\},
$K_{4,2,1,1,1}$-\{$(a_1,b_1),(a_2,c)$\},
$K_{4,2,1,1,1}$-\{$(a_1,b_1),(a_2,b_1)$\},
$K_{4,2,1,1,1}$-\{$(a_1,c),(a_2,c)$\},
$K_{3,3,1,1,1}$-\{$(a_1,b_1),(a_2,b_2)$\},
$K_{3,2,2,1,1}$-\{$(a_1,b_1),(a_1,c_1)$\},
$K_{3,2,2,1,1}$-\{$(a_1,b_1),(a_2,b_2)$\},
$K_{2,2,1,1,1,1}$-\{$(b_1,c),(b_2,d)$\},
$K_{2,2,1,1,1,1}$-\{$(a_1,c),(b_1,c)$\},
$K_{4,1,1,1,1,1}$-\{$(a_1,b),(a_2,b)$\},
$K_{4,1,1,1,1,1}$-\{$(a_1,b),(a_2,c)$\},
$K_{3,2,1,1,1,1}$-\{$(a_1,b_1),(a_2,b_2)$\},
$K_{2,1,1,1,1,1,1}$-\{$(a_1,b),(b,c)$\},
$K_{2,1,1,1,1,1,1}$-\{$(b,c),(d,e)$\}, and
$K_{2,1,1,1,1,1,1}$-\{$(a_1,b),(a_2,c)$\}.

$K_{3,1,1,1,1,1}$-\{$(b,c),(d,e)$\} and
$K_{3,2,1,1,1}$-\{$(b_1,c),(b_2,d)$\} are minors of
$K_{2,1,1,1,1,1,1}$-\{$(b,c),(d,e)$\}.

$K_{4,2,1,1,1}$-\{$(a_1,c),(a_2,d)$\},
$K_{4,3,1,1}$-\{$(a_1,c),(a_2,d)$\},
$K_{4,3,1,1}$-\{$(a_1,b_1),(a_2,c)$\},
$K_{4,2,2,1}$-\{$(a_1,b_1),(a_2,d)$\},
$K_{4,2,2,1}$-\{$(a_1,b_1),(a_2,c_1)$\},
$K_{4,2,2,1}$-\{$(a_1,b_1),(a_2,b_2)$\},
$K_{4,3,2}$-\{$(a_1,c_1),(a_2,c_2)$\}, 
$K_{4,3,2}$-\{$(a_1,b_1),(a_2,c_1)$\}, 
$K_{4,3,2}$-\{$(a_1,b_1),(a_2,b_2)$\}, 
$K_{4,4,1}$-\{$(a_1,b_1),(a_2,c)$\}, and
$K_{4,4,1}$-\{$(a_1,b_1),(a_2,b_2)$\}    are  minors of
$K_{4,1,1,1,1,1}$-\{$(a_1,b),(a_2,c)$\}.

$K_{4,3,1,1}$-\{$(a_1,c),(a_2,c)$\}, $K_{4,2,2,1}$-\{$(a_1,d),(a_2,d)$\},
$K_{4,2,2,1}$-\{$(a_1,b_1),(a_2,b_1)$\}, 
$K_{4,3,2}$-\{$(a_1,c_1),(a_2,c_1)$\},
$K_{4,3,2}$-\{$(a_1,b_1),(a_2,b_1)$\}, and
$K_{4,4,1}$-\{$(a_1,c),(a_2,c)$\} are minors of
$K_{4,1,1,1,1,1}$-\{$(a_1,b),(a_2,b)$\}.

$K_{3,3,2,1}$-\{$(b_1,c_1),(b_2,c_2)$\},
$K_{3,3,2,1}$-\{$(a_1,b_1),(a_2,b_2)$\},
$K_{3,2,2,2}$-\{$(a_1,b_1),(a_2,b_2)$\},
$K_{4,3,2}$-\{$(b_1,c_1),(b_2,c_2)$\} and
$K_{3,3,3}$-\{$(a_1,b_1),(a_2,b_2)$\}  are minors of
$K_{3,2,1,1,1,1}$-\{$(a_1,b_1),(a_2,b_2)$\}.

$K_{3,3,2,1}$-\{$(a_1,b_1),(b_1,c_1)$\},
$K_{3,2,2,2}$-\{$(a_1,b_1),(a_1,c_1)$\},
$K_{4,3,2}$-\{$(a_1,b_1),(a_1,c_1)$\},
$K_{4,3,2}$-\{$(a_1,b_1),(b_1,c_1)$\},
$K_{3,3,3}$-\{$(a_1,b_1),(b_1,c_1)$\}, and
$K_{4,2,2,1}$-\{$(a_1,b_1),(a_1,c_1)$\} are minors of
$K_{3,2,2,1,1}$-\{$(a_1,b_1),(a_1,c_1)$\}.

$K_{4,4,1}$-\{$(a_1,b_1),(a_1,b_2)$\} is not intrinsically knotted by
the corollary to lemma \ref{flemming}. Delete vertices $a_2$ and $c$ for
a planar graph. 
\qed

\subsection{Proof of Adams' conjecture for $2$-deficient graphs}

\begin{thm}\label{Aconj2def} If $G$ is a 2-deficient graph, and any one
vertex and the edges coming into it are removed, the remaining graph is
intrinsically linked.
\end{thm}

\Pf

It suffices to verify the theorem for minimal examples of
2-deficient graphs. 

\begin{description}
\item[$k = 1$] $K_8$-2e is intrinsically knotted; if we remove a vertex,
we get $K_7$-2e, $K_7$-e, or $K_7$, all of which are intrinsically linked.

\item[$k = 2$] $K_{5,5}$-2e is intrinsically knotted; if we remove a
vertex we get $K_{5,4}$-2e, $K_{5,4}$-e, or $K_{5,4}$, all of which are
intrinsically linked.

\item[$k = 3$] If a vertex is removed from a minimal knotted 2-deficient graph that
has 3 parts, the result is one of the following graphs:
$K_{4,3,1}$-\{$(a_1,c),(b_1,c)$\}, $K_{4,3,1}$-\{$(a_1,b_1),(b_1,c)$\},
$K_{4,3,1}$-\{$(b_1,c),(b_2,c)$\}, $K_{4,3,1}$-\{$(a_1,b_1),(b_2,c)$\},
$K_{4,3,1}$-\{$(a_1,b_1),(a_1,b_2)$\},
$K_{4,3,1}$-\{$(a_1,b_1),(a_1,c)$\},
$K_{3,3,2}$-\{$(a_1,b_1),(a_1,b_2)$\},
$K_{3,3,2}$-\{$(a_1,c_1),(a_2,c_1)$\},
$K_{3,3,2}$-\{$(a_1,c_1),(b_1,c_1)$\},
$K_{3,3,2}$-\{$(a_1,b_1),(b_2,c_1)$\},
$K_{3,3,2}$-\{$(a_1,c_1),(b_1,c_2)$\},
$K_{3,3,2}$-\{$(a_1,c_1),(a_1,c_2)$\},
$K_{4,2,2}$-\{$(a_1,b_1),(b_2,c_1)$\},
$K_{4,2,2}$-\{$(a_1,b_1),(b_1,c_1)$\},
$K_{4,2,2}$-\{$(b_1,c_1),(b_2,c_1)$\},
$K_{4,2,2}$-\{$(a_1,b_1),(a_1,b_2)$\}, $K_{5,3,1}$-2e, $K_{4,4,1}$-2e,
$K_{4,3,2}$-2e, $K_{3,3,3}$-2e, $K_{5,2,2}$-2e, $K_{5,4}$-2e,
$K_{4,4}$-e, or one of these graphs with 1 or 2 fewer edges missing. Note
that they are all intrinsically linked.

\item[$k = 4$] If a vertex is removed from a minmal knotted 2-deficient graph that
has 4 parts, the result is one of the following graphs:
$K_{4,3,1}$-\{$(b,c),e$\}, $K_{4,3,1}$-\{$(a_1,b_1),(a_1,b_2)$\},
$K_{4,3,1}$-\{$(a_1,b_1),(a_1,c)$\},
$K_{3,3,2}$-\{$(a_1,b_1),(a_1,b_2)$\},
$K_{3,3,2}$-\{$(a_1,c_1),(a_2,c_1)$\},
$K_{3,3,2}$-\{$(a_1,c_1),(b_1,c_1)$\},
$K_{3,3,2}$-\{$(a_1,b_1),(b_2,c_1)$\},
$K_{3,3,2}$-\{$(a_1,c_1),(b_1,c_2)$\},
$K_{3,3,2}$-\{$(a_1,c_1),(a_1,c_2)$\}, $K_{4,2,2}$-\{$(b,c),e$\},
$K_{4,2,2}$-\{$(a_1,b_1),(a_1,b_2)$\}, $K_{5,3,1}$-2e, $K_{4,4,1}$-2e,
$K_{4,3,2}$-2e, $K_{3,3,3}$-2e, $K_{5,2,2}$-2e,
$K_{4,2,1,1}$-\{$(b,c),e$\}, $K_{4,2,1,1}$-\{$(c,d),e$\},
$K_{4,2,1,1}$-\{$(a_1,b_1),(a_1,b_2)$\},
$K_{4,2,1,1}$-\{$(a_1,b_1),(a_1,c)$\},
$K_{4,2,1,1}$-\{$(a_1,c),(a_1,d)$\}, $K_{3,3,1,1}$-\{$(b,c),e$\},
$K_{3,3,1,1}$-\{$(c,d),e$\}, $K_{3,3,1,1}$-\{$(a_1,b_1),(a_1,b_2)$\},
$K_{3,2,2,1}$-\{$(b,c),e$\}, $K_{3,2,2,1}$-\{$(c,d),e$\},
$K_{3,2,2,1}$-\{$(a,d),e$\}, $K_{3,2,2,1}$-\{$(a_1,b_1),(a_1,b_2)$\},
$K_{3,2,2,1}$-\{$(a_1,b_1),(a_2,b_1)$\},
$K_{3,2,2,1}$-\{$(a_1,b_1),(a_2,c_1)$\}, $K_{2,2,2,2}$-2e,
$K_{5,2,1,1}$-2e, $K_{4,3,1,1}$-2e, $K_{4,2,2,1}$-2e, $K_{3,3,2,1}$-2e,
$K_{3,2,2,2}$-2e, or one of these graphs with 1 or 2 fewer edges missing.
Note that they are all intrinsically linked.

\item[$k = 5$] If a vertex is removed from a minimal knotted 2-deficient graph that
has 5 parts, the result is one of the following graphs:
$K_{3,2,1,1}$-\{$(b_1,c),(b_1,d)$\}, $K_{3,2,1,1}$-\{$(b_1,c),(b_2,c)$\},
$K_{4,2,1,1}$-\{$(b,c),e$\}, $K_{4,2,1,1}$-\{$(c,d),e$\},
$K_{4,2,1,1}$-\{$(a_1,b_1),(a_1,b_2)$\},
$K_{4,2,1,1}$-\{$(a_1,b_1),(a_1,c)$\},
$K_{4,2,1,1}$-\{$(a_1,c),(a_1,d)$\}, $K_{3,3,1,1}$-\{$(b,c),e$\},
$K_{3,3,1,1}$-\{$(c,d),e$\}, $K_{3,3,1,1}$-\{$(a_1,b_1),(a_1,b_2)$\},
$K_{3,2,2,1}$-\{$(b,c),e$\}, $K_{3,2,2,1}$-\{$(c,d),e$\},
$K_{3,2,2,1}$-\{$(a,d),e$\}, $K_{3,2,2,1}$-\{$(a_1,b_1),(a_1,b_2)$\},
$K_{3,2,2,1}$-\{$(a_1,b_1),(a_2,b_1)$\},
$K_{3,2,2,1}$-\{$(a_1,b_1),(a_2,c_1)$\}, $K_{2,2,2,2}$-2e,
$K_{5,2,1,1}$-2e, $K_{4,3,1,1}$-2e, $K_{4,2,2,1}$-2e, $K_{3,3,2,1}$-2e,
$K_{2,2,1,1,1}$-\{$(b_1,c),(b_2,c)$\},
$K_{2,2,1,1,1}$-\{$(b_1,c),(b_1,d)$\}, $K_{3,1,1,1,1}$-\{$(b,c),(c,d)$\},
$K_{4,1,1,1,1}$-\{$(b,c),e$\}, $K_{4,1,1,1,1}$-\{$(a_1,b),(a_1,c)$\},
$K_{3,2,1,1,1}$-\{$(a,c),e$\}, $K_{3,2,1,1,1}$-\{$(b,c),e$\},
$K_{3,2,1,1,1}$-\{$(c,d),e$\}, $K_{3,2,1,1,1}$-\{$(a_1,b_1),(a_1,b_2)$\},
$K_{3,2,1,1,1}$-\{$(a_1,b_1),(a_2,b_1)$\}, $K_{2,2,2,1,1}$-2e,
$K_{5,1,1,1,1}$-2e, $K_{4,2,1,1,1}$-2e, $K_{3,3,1,1,1}$-2e,
$K_{3,2,2,1,1}$-2e, or one of these graphs with 1 or 2 fewer edges missing.
Note that they are all intrinsically linked.

\item[$k = 6$] If a vertex is removed from a minimal knotted 2-deficient graph that
has 6 parts, the result is one of the following graphs:
$K_{2,2,1,1,1}$-\{$(b_1,c),(b_2,c)$\},
$K_{2,2,1,1,1}$-\{$(a_1,c),(b_1,d)$\},
$K_{2,2,1,1,1}$-\{$(b_1,c),(b_1,d)$\}, $K_{3,1,1,1,1}$-\{$(b,c),(c,d)$\},
$K_{4,1,1,1,1}$-\{$(b,c),e$\}, $K_{4,1,1,1,1}$-\{$(a_1,b),(a_1,c)$\},
$K_{3,2,1,1,1}$-\{$(a,c),e$\}, $K_{3,2,1,1,1}$-\{$(b,c),e$\},
$K_{3,2,1,1,1}$-\{$(c,d),e$\}, $K_{3,2,1,1,1}$-\{$(a_1,b_1),(a_1,b_2)$\},
$K_{3,2,1,1,1}$-\{$(a_1,b_1),(a_2,b_1)$\}, $K_{2,2,2,1,1}$-2e,
$K_{5,1,1,1,1}$-2e, $K_{4,2,1,1,1}$-2e, $K_{3,3,1,1,1}$-2e,
$K_{2,1,1,1,1,1}$-\{$(a_1,b),(a_1,c)$\},
$K_{2,1,1,1,1,1}$-\{$(a_1,b),(a_2,b)$\},
$K_{2,1,1,1,1,1}$-\{$(a_1,b),(c,d)$\},
$K_{2,1,1,1,1,1}$-\{$(b,c),(c,d)$\}, $K_{3,1,1,1,1,1}$-2e,
$K_{2,2,1,1,1,1}$-2e, $K_{4,1,1,1,1,1}$-2e, $K_{3,2,1,1,1,1}$-2e, or one
of these graphs with 1 or 2 fewer edges missing. Note that they are all
intrinsically linked.

\item[$k = 7$] If a vertex is removed from a minimal knotted 2-deficient graph that
has 7 parts, the result is one of the following graphs:
$K_{2,1,1,1,1,1}$-\{$(a_1,b),(a_1,c)$\},
$K_{2,1,1,1,1,1}$-\{$(a_1,b),(a_2,b)$\},
$K_{2,1,1,1,1,1}$-\{$(a_1,b),(c,d)$\},
$K_{2,1,1,1,1,1}$-\{$(b,c),(c,d)$\}, $K_{3,1,1,1,1,1}$-2e,
$K_{2,2,1,1,1,1}$-2e, $K_{1,1,1,1,1,1,1}$-2e, $K_{2,1,1,1,1,1,1}$-2e,  or
one of these graphs with 1 or 2 fewer edges missing. Note that they are
all intrinsically linked.

\item[$k \geq 8$] If a vertex is removed from a knotted 2-deficient graph
that has 8 or more parts, the result is a 2-deficient, 1-deficient, or
complete partite graph with 7 or more parts, all of which are
intrinsically linked.

\end{description}

\qed

\section{Graphs on $8$ vertices}
In this section we provide a classification of intrinsically knotted
graphs on $8$ vertices. We verify that removing a vertex from
any of these results in an intrinsically linked graph. We also discuss
a question of Sachs~\cite{S} about maximal unlinked and unknotted graphs.

\subsection{Classification}

Graphs on $8$ vertices are subgraphs of $K_8$. We will examine in turn
subgraphs which are obtained by removing $1$, $2$, $3$, $\ldots$ edges 
from $K_8$.  We have already noted that
$K_8$, $K_8 - e$, and both $K_8-2e$ graphs are intrinsically
knotted. Of the five graphs $K_8-3e$, the three which are intrinsically 
knotted can all be obtained by removing two edges from
$K_{2,1,1,1,1,1,1}$. (See classification of $2$-deficient graphs
above. Note that $2,1,1,1,1,1,1-\{(a,b),(c,d)\}$ and
$2,1,1,1,1,1,1 - \{(b,c),(c,d) \}$ are the same graph.)

There are $11$ graphs of the form $K_8-4e$. Seven of these are not
knotted as they can be realized by removing an edge from one of the two
unknotted
$K_8-3e$ graphs. The remaining four are intrinsically knotted. Three of
these four are  of the form $K_{2,2,1,1,1,1} - 2e$. The fourth is
obtained by removing
$4$ edges all incident to the same vertex. This graph is intrinsically
knotted as it contains $K_7$ as a minor.

There are $24$ graphs $K_8-5e$. All but $4$ of these are  not
knotted as they are minors of an unknotted $K_8-4e$. The four
intrinsically knotted $K_8-5e$'s are perhaps most easily described in
terms of their complementary graphs.
\begin{figure}[ht]
\begin{center}
\includegraphics[scale=0.5]{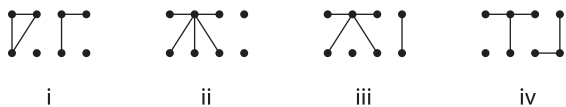}

\caption{Intrinsically knotted $K_8-5e$'s}\label{fig5def}
\end{center}
\end{figure} 
In figure~\ref{fig5def}, Graph i is
$K_{3,2,1,1,1} - (b,c)$, or, equivalently, $K_{3,1,1,1,1,1} - \{(b,c),(c,d)\}$.
The other three graphs in the figure have $K_7$ as a minor and are,
therefore, intrinsically knotted.

Of the $56$ $K_8 - 6e$ graphs, all but $6$ are minors of unknotted
$K_8 - 5e$'s. 
\begin{figure}[ht]
\begin{center}
\includegraphics[scale=0.15]{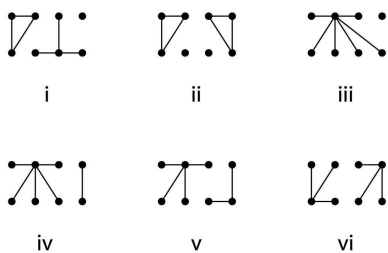}

\caption{Intrinsically knotted $K_8-6e$'s}\label{fig6def}
\end{center}
\end{figure} 
In figure~\ref{fig6def}, Graph i is $K_{3,2,1,1,1} -
\{(b_1,c),(b_1,d)\}$ while Graph ii is $K_{3,3,1,1}$, or, equivalently,
$K_{3,2,1,1,1} - \{(b_1,c),(b_2,c)\}$. The remaining four graphs are
obtained by splitting a vertex of $K_7$  and are therefore intrinsically
knotted.

Only $2$ of the $K_8-7e$ graphs are not minors of some
unknotted $K_8 -6e$. One of these two is $H_8$ \cite{KS}. The other is
$K_7$ with one additional vertex. These are both intrinsically knotted. 
Moreover, $K_7$ and $H_8$ are minor minimal \cite{KS}. Thus any subgraph of $K_8$
obtained by removing $8$ or more edges is not intrinsically knotted. 

In total then, there are twenty intrinsically knotted graphs on $8$
vertices. 

\subsection{Proof of Adams' conjecture for graphs on $8$ vertices}

\begin{thm} If $G$ is an intrinsically knotted graph on $8$ vertices and
any one vertex and the edges coming into it are removed, the remaining 
graph is intrinsically linked.
\end{thm}

\Pf
We have already verified this for the knotted graphs that are $0$-, $1$-,
or $2$-deficient complete partite graphs. Most of the other knotted graphs
have $K_7$ as a minor. On removing a vertex, the resulting graph will have
$K_6$ as a minor and be intrinsically linked. 
Foisy has shown that the removal of a vertex from
$H_8$ results in an intrinsically linked graph \cite{F2}.



\qed

\subsection{A question of Sachs}

Sachs~\cite{S} asked if a graph on $n$ vertices that is not intrinsically
linked could have more than $4n-10$ edges. Using lemma~\ref{flemming}, 
we can ask a 
similar question about intrinsically knotted graphs. For $n \geq 5$, a 
planar triangulation with $n-2$ vertices will have $3(n-4)$ edges. Adding
a $K_2$ gives a graph with $5n-15$ edges that is not intrinsically 
knotted by lemma~\ref{flemming}. This is the maximum for $n = 5,6,7$, 
and our analysis
of graphs on $8$ vertices shows that it is also the maximum for $n=8$. 
In other words, a graph with more than $5n-15$ edges on $n$ vertices
is intrinsically knotted when $5 \leq n \leq 8$.
Is this also true for larger $n$?

\noindent%
{\bf Question:} Is there a graph on $n$ vertices that is not
intrinsically knotted and has more than $5n-15$ edges?


\end{document}